\documentclass[12pt]{amsart}
\usepackage{amssymb,eucal,graphicx}

\usepackage[all,poly,import,color,ps,dvips]{xy}

\SelectTips{cm}{}

\numberwithin{equation}{section}

\newtheorem{theorem}[equation]{Theorem}
\newtheorem*{theorema}{Theorem}
\newtheorem{coroll}[equation]{Corollary}
\newtheorem{corollary}[equation]{Corollary}
\newtheorem*{corollario}{Corollary}

\newtheorem{proposition}[equation]{Proposition}

\theoremstyle{definition}
\newtheorem{definition}[equation]{Definition}
\newtheorem{notation}[equation]{Notation}
\newtheorem{example}[equation]{Example}
\newtheorem*{question}{Question}

\newtheorem{remark}[equation]{Remark}

\DeclareMathOperator{\coker}{coker}

\DeclareMathOperator{\Sing}{Sing} \DeclareMathOperator{\im}{Im}
\DeclareMathOperator{\tr}{tr} \DeclareMathOperator{\sgn}{sgn}
\DeclareMathOperator{\Sym}{Sym}

\newcommand{\mapright}[1]{\stackrel{#1}\longrightarrow}

\newcommand\Oc{{\mathcal O}}

\newcommand\X{\mathcal X}
\newcommand\Ha{\mathcal H}
\newcommand\D{\Delta}

\newcommand\Z{\mathcal Z}

\newcommand\CC{\mathbb{C}}
\newcommand\NN{\mathbb{N}}

\newcommand\R{\mathcal R}
\newcommand\T{\mathcal T}

\newcommand\Pp{\mathbb P}
\newcommand\PR{\mathbb{P}^r}

\renewcommand{\P}{\mathbb{P}}
\renewcommand{\geq}{\geqslant}
\renewcommand{\leq}{\leqslant}
\renewcommand{\ge}{\geqslant}
\renewcommand{\le}{\leqslant}

\newcommand{\st}{\scriptstyle}
\newcommand{\dt}{\displaystyle}

\def\colorxy(#1){/xycolor{#1 setrgbcolor}def}

\title[On the geometric genus of reducible surfaces]{On the geometric genus 
of reducible surfaces \\
and degenerations of surfaces to unions of planes}

\author{A.~Calabri, C.~Ciliberto, F.~Flamini, R.~Miranda}

\email{calabri@mat.uniroma2.it} \curraddr{Dipartimento di
Matematica, Universit\`a degli Studi di Roma ``Tor Vergata",
Via della Ricerca Scientifica, 00133 Roma, Italy}

\email{cilibert@mat.uniroma2.it} \curraddr{Dipartimento di
Matematica, Universit\`a degli Studi di Roma ``Tor Vergata",
Via della Ricerca Scientifica, 00133 Roma, Italy}

\email{flamini@mat.uniroma3.it} \curraddr{Dipartimento di
Matematica Pura ed Applicata, Universit\`a degli Studi di L'Aquila,
Via Vetoio, Loc.~Coppito, 67100 L'Aquila, Italy}

\email{miranda@math.colostate.edu} \curraddr{Department of
Mathematics, 101 Weber Building, Colorado State University,
Fort Collins, CO 80523--1874, U.S.A.}

\thanks{
The first two authors are partially supported by E.C.\ project
EAGER, contract n.\ HPRN-CT-2000-00099.}

\begin{document}

\maketitle

\section{Introduction}\label{S:1}

In this paper we study some properties of degenerations
of surfaces whose general fibre is a smooth projective algebraic surface
and whose central fibre is a reduced, connected surface $X \subset \PR$, $r \geq 3$, 
which is assumed to be a union of planes. Here 
we present a first set of results on the subject; other aspects are still work in 
progress and will appear later (see \cite{CCFM}).  

Our original motivation has been a series of papers by Guido Zappa which appeared in the 1940--50's 
regarding degenerations of scrolls to unions of planes and the computation of bounds
for the topological invariants of an arbitrary smooth projective surface
which is assumed to degenerate to a union of planes (see \cite{Za1, Za1b,
Za1c, Za2, Za2b, Za3b, Za3} and \cite{CCFM}).

Zappa was in turn motivated by earlier papers by Francesco
Severi concerning \emph{Zeuthen's problem}, i.e.\ the existence of degenerations
of smooth projective (space) curves to unions of lines with only nodes as singularities
(now called \emph{stick curves}).

Zeuthen's problem has been studied by several authors, also recently
(see e.g.\ \cite{Hart});
on the contrary, unions of planes have been studied only in terms of degenerations of a few types
of smooth surfaces, e.g.\ K3 surfaces (see \cite{CLM, CMT, CMT2, FM}).

In this paper, we first
study the geometry and the combinatorics of a union of planes $X$
considered as a reduced, connected surface on its own (cf. \S \ref{S:3} and \ref{S:4pg}); then,
we focus on the
case in which $X$ is the central fibre of an embedded
degeneration $\X\to\D$, where $\D$ is the complex unit
disk and where $\X \subseteq \D \times \PR$, $r \geq 3$, is a closed subscheme
of relative dimension two. In this case, we deduce some
properties of the general fibre $\X_t$, $t \neq 0$, of the degeneration from the
ones of its central fibre $\X_0 = X$ (see \S \ref{S:zapdeg}).

It is well-known that, in dimension one, for any integer $g\ge2$ any smooth projective curve
of genus $g$ with general moduli and sufficiently general degree
can be degenerated to a suitable stick curve (see, e.g. \cite{AW} and \cite{Sev}).

On the contrary, in dimension two, worse singularities than normal crossings are needed
in order to degenerate as many surfaces
as possible to unions of planes (cf.\ \cite{CCFM}).

Here we shall focus on the case of $X$ a union of planes --- or
more generally a union of smooth projective surfaces --- whose
singularities are:
\begin{itemize}
\item in codimension one, double curves which are smooth and irreducible;
\item multiple points, which are locally analytically isomorphic
to the vertex of a cone over a stick curve
with arithmetic genus either zero or one and which is projectively
normal in the projective space it spans.
\end{itemize}
These multiple points will be called \emph{Zappatic singularities},
whereas a surface like $X$ will be called a \emph{Zappatic surface}.
If moreover $X \subset \PR$, for some positive $r$, and if all its irreducible components
are planes, then $X$ is said to be a \emph{planar Zappatic surface}.

Actually we will concentrate on the so called \emph{good Zappatic surfaces},
i.e.\ Zappatic surfaces having only Zappatic singularities whose associated
stick curve has one of the following
dual graphs (cf.\ Examples \ref{ex:tngraphs} and \ref{ex:zngraphs}, Definition \ref{def:goodzapp},
Figures \ref{fig:zapsings} and \ref{fig:R3E3R4}):
\begin{itemize}
    \item[$R_n$:] a chain of length $n$, with $n\ge3$;
    \item[$S_n$:] a fork with $n-1$ teeth, with $n\ge4$;
    \item[$E_n$:] a cycle of order $n$, with $n\ge3$.
\end{itemize}
Let us call $R_n$-, $S_n$-, $E_n$-{\em point} the corresponding
multiple point of the Zappatic surface $X$. These singularities 
play a major role in the whole subject (cf. \cite{CCFM}).

We remark that a Zappatic surface $X$ is locally Gorenstein
(i.e.\ its dualizing sheaf $\omega_X$ is invertible)
if and only if its Zappatic singularities are only $E_n$-points, for any $n\ge3$.

We associate to a good Zappatic surface $X$ a graph $G_X$
(see Definition \ref{def:dualgraph})
which encodes the configuration of the irreducible components of $X$
as well as its Zappatic singularities.

We shall see (cf.\ Sections \ref{S:3} and \ref{S:4pg}) how to combinatorially compute 
from the associated
graph $G_X$ some intrinsic and extrinsic invariants of $X$,
e.g.\ the Euler-Poincar{\'e} characteristic $\chi(\Oc_X)$,
the geometric genus $p_g(X)$
(cf.\ Remark \ref{rem:pg}), as well as
--- when $X \subset \PR$, $r \geq 3$ --- the degree
$d=\deg(X)$, the sectional genus $g$, and so on.

When $X$ is further assumed to be the central fibre of a degeneration $\X\to\D$
(resp., of an embedded degeneration $\X \to \D$, where
$\X \subset \D \times \PR$, $r \geq 3$)
we will then compute intrinsic (resp., intrinsic and extrinsic) invariants
of the general fibre $\X_t$, for $t \neq 0$.

We shall see how to directly compute some of the invariants of $X$ by means of the associated
graph $G_X$. Determining formulas for a few invariants (e.g.\ $d$ and $g$) is quite easy,
whereas for other invariants, like $\chi(\Oc_X)$, it requires some more work.

Actually the computation of the geometric genus is still an open question, in general.
Indeed, we first prove the following (cf.\ Theorem \ref{thm:4.pgbis}):

\begin{theorema}
Let $X=\bigcup_i X_i$ be a Zappatic surface with global normal
crossings, i.e.\ with only $E_3$-points as Zappatic singularities.
Denote by $\omega_X$ the dualizing sheaf of $X$ and by $G_{X}$ the
associated graph of $X$.
Consider the natural map:
$$\Phi:\bigoplus_i H^1(X_i,\Oc_{X_i}) \to \bigoplus_{i,j} H^1(C_{ij},\Oc_{C_{ij}}),$$
where $C_{ij}=X_i\cap X_j$ (cf.\ formula \eqref{eq:4.fi}).
Then, the following inequality holds:
\begin{equation}\label{eq:intro}
p_g(X):=h^0(X,\omega_X) \leq h_2(G_X,\CC)+\sum_{i=1}^v p_g(X_i)+ \dim(\coker(\Phi)).
\end{equation}

Furthermore, a sufficient condition for the equality in \eqref{eq:intro}
to hold is either that 
\begin{itemize}
\item[(i)] each irreducible component $X_i$ is a regular surface (i.e.\ $h^1(\Oc_{X_i})=0$), or that 
\item[(ii)] for any irregular component $X_j$ of $X$,
the divisor $C_{j} := X_{j} \cap \overline{(X \setminus X_{j})}$ is ample on $X_j$.
\end{itemize}
\end{theorema}

The proof of the above theorem also shows the following:

\begin{corollario}
Let $X$ be a planar Zappatic surface with global normal crossings (i.e.\ only $E_3$-points)
and $G_X$ be its associated graph.
Then, there exists an explicit isomorphism
\begin{equation}\label{eq:intro2}
H^0(X,\omega_X) \cong H_2(G_X,\CC),
\end{equation}
where $\omega_X$ is the dualizing sheaf of $X$.
Therefore
\[
p_g(X):= h^0(X,\omega_X) = h_2(G_X,\CC).
\]
\end{corollario}

From the proof of Theorem \ref{thm:4.pgbis}, it will be clear that
the isomorphism \eqref{eq:intro2} essentially follows from
evaluation of residues at the $E_3$-points, with a suitable use of signs.

We show that equality holds in \eqref{eq:intro} when $X$
is smoothable,
i.e.\ when $X$ is the central fibre of a {\em semistable degeneration} $\X\to\D$:
this follows by the computation of
the geometric genus $p_g(\X_t)$ of the
general fibre $\X_t$, for $t \neq 0$, via the Clemens-Schmid exact sequence
(cf.\ Theorem \ref{thm:4.CSbis}).

We remark that our computation of the geometric genus is independent of the fact that
$X$ is the central fibre of a semistable
degeneration. We deal with this particular case in \S
\ref{S:zapdeg}, where we show that in a semistable 
degeneration, whose central fibre is a
Zappatic surface with only $E_3$-points, the geometric genus of the fibres is constant
(see Corollary \ref{cor:4.CSganzo}).

It is still an open problem to find an example 
for which the strict inequality holds in \eqref{eq:intro}.

Finally, we will see that the above results can be generalized to a smoothable
good Zappatic surface, i.e.\ with $R_n$-, $S_n$- and $E_n$-points, for any $n\ge3$
(see Theorem \ref{thm:pg}).

A natural question to ask is which Zappatic singularities
are needed in order to degenerate as many surfaces as possible.
Results and some examples contained in \cite{CCFM} suggest that,
even if a given projective surface $X$ needs $E_n$-, $R_n$-, or $S_n$-points with large $n$,
there might be a birational model of $X$ which needs just $R_3$- and $E_n$-points, with 
$n\le 6$. For example, in \cite{CMT} there are interesting examples of K3 surfaces degenerating to
a Zappatic surface with at most $R_3$- and $E_6$-points, called \emph{pillow} degenerations.
However, we do not have enough evidence to state a reasonable conjecture in this direction.

\bigskip
{\it Acknowledgments.} The authors would like to thank Janos Koll\'ar, for some useful discussions and
references, and the organizers of the Fano Conference,
for the very stimulating atmosphere during the whole week of the meeting.

\bigskip

\section {Reducible curves and associated graphs}\label{S:2}

Let $C$ be a projective curve and let $C_i$, $i=1,\ldots,n,$ be
its irreducible components. We will assume that:

\begin{itemize}
\item $C$ is connected and reduced;
\item $C$ has at most nodes as singularities;
\item the curves $C_i, i=1,\ldots,n,$ are smooth.
\end{itemize}

If two components $C_i, \; C_j, \; i<j,$ intersect at $m_{ij}$
points, we will denote by $P_{ij}^h, \; h=1,\ldots,m_{ij}$, the
corresponding nodes of $C$.

We can associate to this situation a simple (i.e.\ with no loops),
connected graph $G_C$:
\begin{itemize}

\item whose vertices  $v_1,\ldots,v_n,$ correspond to the
components $C_1$, $\ldots$, $C_n$;

\item whose edges $\eta^h_{ij}$, $i<j, \; h=1,\ldots,m_{ij}$,
joining the vertices $v_i$ and $v_j$, correspond to the nodes
$P_{ij}^h$ of $C$.
\end{itemize}

We will assume the graph to be {\em lexicographically oriented},
i.e.\ each edge is assumed to be oriented from the vertex with
lower index to the one with higher index.

We will use the following notation:

\begin{itemize}
\item $v:$ the number of vertices of $G_C$, i.e.\ $v=n$;
\item $e:$ the number of edges of $G_C$, i.e.\ the number of nodes of $C$;
\item $g_i:$ the  genus of the curve $C_i$,
which we consider as the {\em weight} of the vertex $v_i$;
\item $\chi(G_C)=v-e$ is the Euler-Poincar\'e characteristic of $G_C$;
\item $h_1(G_C)=1-\chi(G_C)$ is the first Betti number of $G_C$.
\end{itemize}

Remark that conversely, given any simple, connected, weighted
(oriented) graph $G$, there is some curve $C$ such that $G=G_C$.

One has the following basic result:

\begin{theorem}\label{thm:pacurves} In the above situation
\begin{equation}\label{eq:chicurves}
\chi(\Oc_C)=\chi(G_C) - \sum_{i=1}^v g_i=v-e- \sum_{i=1}^v g_i.
\end{equation}
\end {theorem}

\begin{proof} Let $\nu: \tilde{C} \to C$ be the normalization morphism; this defines
the exact sequence of sheaves on $C$:
\begin{equation}\label{eq:4gen}
0 \to \Oc_C \to \nu_* (\Oc_{\tilde{C}}) \to \underline{\tau} \to 0,
\end{equation}
where $\underline{\tau} $ is a sky-scraper sheaf supported at $\Sing(C)$.
Since the singularities of $C$ are only nodes, one easily determines
$H^0(C, \underline{\tau} ) \cong \CC^e$. Therefore, by the exact sequence \eqref{eq:4gen}, one gets
$$\chi(\Oc_C) = \chi(\nu_* (\Oc_{\tilde{C}})) -e.$$By the Leray isomorphism and by the
fact that $\nu$ is finite, one has
$\chi(\nu_* (\Oc_{\tilde{C}})) = \chi(\Oc_{\tilde{C}})$. Since $\tilde{C}$ is
a disjoint union of the $v=n$ irreducible components of $C$, one has
$\chi(\Oc_{\tilde{C}}) = v - \sum_{i=1}^v g_i$, which proves \eqref{eq:chicurves}.
(Cf.\ also \cite{Bar} for another proof.)
\end{proof}

We remark that formula \eqref{eq:chicurves} is equivalent to
\begin{equation} \label{eq:genuscurves}
p_a(C)= h_1(G_C)+ \sum_{i=1}^v g_i,
\end{equation}(cf.\ Proposition \ref{prop:ghypsect}).

Notice that $C$ is locally Gorenstein, i.e.\ the dualizing
sheaf $\omega_C$ is invertible.
One defines the {\it geometric genus} of $C$ to be
\begin{equation}\label{eq:ggcurves}
p_g(C):= h^0(C,\omega_C).
\end{equation}

By the Riemann-Roch Theorem, one has
\begin{equation}
p_g(C)=p_a(C)=h_1(G_C)+\sum_{i=1}^v g_i=e-v+1+\sum_{i=1}^v g_i.
\end{equation}

However, one can prove the previous formula combinatorially by
showing that there is a natural short exact sequence:
\begin{equation}
0\to \bigoplus_{i=1}^vH^0(C_i, \omega_{C_i}) \to H^0(C,\omega_C) \to
H_1(G_C,\CC)\to 0
\end{equation}

We will not dwell on this now, since we shall show the existence
of an analogous sequence in the surface case in \S \ref {S:4pg}.

If we have a flat family ${\mathcal C}\to \D$ over a disc $\D$
with general fibre ${\mathcal C}_t$ a smooth and irreducible curve
of genus $g$ and special fibre ${\mathcal C}_0=C$, then we can
combinatorially compute $g$ via the formula:
$$
g=p_a(C)=h_1(G_C)+\sum_{i=1}^vg_i=e-v+1+\sum_{i=1}^v g_i.
$$

Usually we will consider a curve $C$ embedded in a projective space
$\PR$. In this situation each curve $C_i$ will have a certain
degree $d_i$, and we will consider the graph $G_C$ as {\em double
weighted}, by attaching to each vertex the pair of weights
$(g_i,d_i)$. Moreover we will attribute to the graph a further
marking number, i.e.\ the embedding dimension $r$ of $C$.

The total degree of $C$ is
$$
d=\sum_{i=1}^v d_i
$$
which is also invariant by flat degeneration.

If each curve $C_i$ is a line, the curve $C$ is called a {\it stick curve}.
In this case the double weighting is $(0,1)$ for each vertex, and
it will be omitted if no confusion arises.

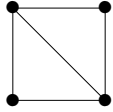
\begin{figure}[ht]
\[
\begin{xy}
0; <35pt,0pt>: 
(0,1)*=0{\bullet};(0,0)*=0{\bullet}**@{-};
(1,0)*=0{\bullet}**@{-}="a";(1,1)*=0{\bullet}**@{-};(0,1)**@{-};"a"**@{-} ,
\end{xy}
\]
\caption{Dual graph of an ``impossible" stick curve.}\label{fig:nostick}
\end{figure}

It should be stressed  that it is not true that for any simple,
connected, double weighted graph $G$  there is a curve $C$ in a
projective space such that $G_C=G$. For example there is no stick
curve corresponding to the graph of Figure \ref{fig:nostick}.

We now give two examples of stick curves which will be frequently used in this paper.

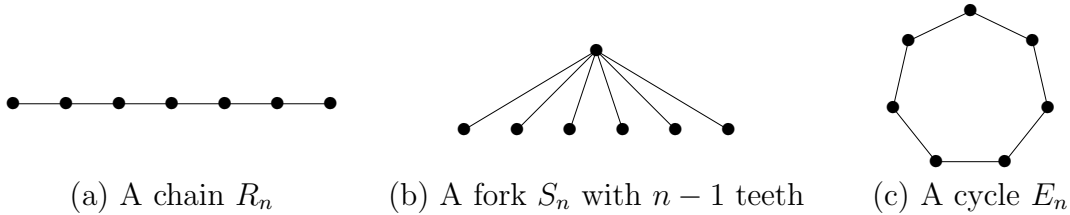
\begin{figure}[h]
\[
\begin{array}{ccc}
\quad
\raisebox{10pt}{$ %
\begin{xy}
0; <20pt,0pt>: 
(0,0)*=0{\bullet};(1,0)*=0{\bullet}**@{-};(2,0)*=0{\bullet}**@{-};(3,0)*=0{\bullet}**@{-};
(4,0)*=0{\bullet}**@{-};(5,0)*=0{\bullet}**@{-};(6,0)*=0{\bullet}**@{-}
,
\end{xy}
$}
\quad & \qquad
\begin{xy}
0; <20pt,0pt>: 
(0,1.5)*=0{\bullet}="u" ,
(-2.5,0)*=0{\bullet};"u"**@{-} ,
(-1.5,0)*=0{\bullet};"u"**@{-} ,
(-0.5,0)*=0{\bullet};"u"**@{-} ,
(2.5,0)*=0{\bullet};"u"**@{-} ,
(1.5,0)*=0{\bullet};"u"**@{-} ,
(0.5,0)*=0{\bullet};"u"**@{-} ,
\end{xy}
\qquad & \qquad
\raisebox{15pt}{$ %
\begin{xy}
0; <30pt,0pt>: 
{\xypolygon7{\bullet}}
\end{xy}
$}
\qquad
\\[4mm]
\text{(a) A chain $R_n$} & \text{(b) A fork $S_n$ with $n-1$ teeth} & \text{(c) A cycle $E_n$}
\end{array}
\]
\caption{Examples of dual graphs.}\label{fig:zapsings}
\end{figure}

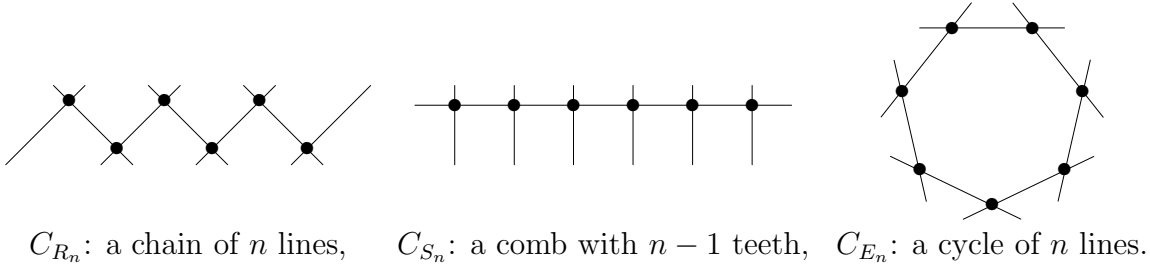
\begin{figure}[h]
\[
\begin{array}{ccc}
\begin{xy}
0; <15pt,0pt>: 
(1.6,1.6)*=0{\bullet},
(2.8,0.4)*=0{\bullet},
(4.0,1.6)*=0{\bullet},
(5.2,0.4)*=0{\bullet},
(6.4,1.6)*=0{\bullet},
(7.6,0.4)*=0{\bullet},
(0,0);(2,2)**@{-} ,
(1.2,2);(3.2,0)**@{-} ,
(2.4,0);(4.4,2)**@{-} ,
(3.6,2);(5.6,0)**@{-} ,
(4.8,0);(6.8,2)**@{-} ,
(6.0,2);(8.0,0)**@{-} ,
(7.2,0);(9.2,2)**@{-} ,
\end{xy}
&
\begin{xy}
0; <15pt,0pt>: 
(0,1.5)*=0{\bullet},
(1.5,1.5)*=0{\bullet},
(3,1.5)*=0{\bullet},
(4.5,1.5)*=0{\bullet},
(6,1.5)*=0{\bullet},
(7.5,1.5)*=0{\bullet},
(-1,1.5);(8.5,1.5)**@{-} ,
(0,0);(0,2)**@{-} ,
(1.5,0);(1.5,2)**@{-} ,
(3,0);(3,2)**@{-} ,
(4.5,0);(4.5,2)**@{-} ,
(6.0,0);(6.0,2)**@{-} ,
(7.5,0);(7.5,2)**@{-} ,
\end{xy}
&
\raisebox{20pt}{$ %
\begin{xy}
0; <35pt,0pt>: 
0,
{\xypolygon7"B"{~:{(-1,0):}~*{\bullet}~>{}}} ,
{\xypolygon7"D"{~={285}~:{(1.2,0):}~>{}}} ,
{\xypolygon7"E"{~={255}~:{(1.2,0):}~>{}}} ,
"E7";"D1"**@{-} ,
"E6";"D7"**@{-} ,
"E5";"D6"**@{-} ,
"E4";"D5"**@{-} ,
"E3";"D4"**@{-} ,
"E2";"D3"**@{-} ,
"E1";"D2"**@{-} ,
\end{xy}
$}
\\[6mm]
\text{$C_{R_n}$: a chain of $n$ lines,}
 & \text{$C_{S_n}$: a comb with $n-1$ teeth,}
 & \text{$C_{E_n}$: a cycle of $n$ lines.}
\end{array}
\]
\caption{Examples of stick curves.}\label{fig:zapsings2}
\end{figure}

\begin{example}\label{ex:tngraphs}
Let $T_n$ be any connected tree with $n \geq 3$ vertices.
This corresponds to a non-degenerate stick curve of degree $n$ in
$\P^n$, which we denote by $C_{T_n}$.
Indeed one can check that, taking a general point $p_i$
on each component of $C_{T_n}$, the line bundle $\Oc_{C_{T_n}} (p_1+\cdots+p_n)$
is very ample.
Of course $C_{T_n}$ has arithmetic genus $0$ and is a
flat limit of rational normal curves in $\P^n$.

We will often consider two particular trees $T_n$:
a chain $R_n$ of length $n$ and the fork $S_n$ with $n-1$ teeth,
i.e.\ a tree consisting of $n-1$ vertices
joining a further vertex (see Figures \ref{fig:zapsings}.(a) and (b)).
The curve $C_{R_n}$ is the union of $n$ lines $l_1,
l_2, \ldots, l_n$ spanning $\P^n$, such that $l_i\cap l_j=\emptyset$ if and only if $1<|i-j|$.
The curve $C_{S_n}$ is the union of $n$
lines $l_1, l_2, \ldots, l_n$ spanning $\P^n$, such that
$l_1,\ldots,l_{n-1}$ all intersect $l_n$ at distinct points (see Figure \ref{fig:zapsings2}).
\end{example}

\vspace{-.6cm}

\begin{example}\label{ex:zngraphs}
Let $Z_n$ be any  simple, connected graph  with $n\geq 3$ vertices
and $h^1(Z_n, \CC)=1$. This corresponds to a projectively normal stick
curve of degree $n$ in $\P^{n-1}$, which we denote by $C_{Z_n}$
(as in Example \ref{ex:tngraphs}).
The curve $C_{Z_n}$ has arithmetic genus $1$ and it is a flat
limit of elliptic normal curves in $\P^{n-1}$.

We will often consider the particular case of a cycle $E_n$ of order $n$ (see Figure \ref{fig:zapsings}.c).
The curve $C_{E_n}$ is the union of $n$ lines $l_1, l_2, \ldots, l_n$
spanning $\P^{n-1}$, such that $l_i\cap l_j=\emptyset$ if and only if $1<|i-j|<n-1$
(see Figure \ref{fig:zapsings2}).

We remark that $C_{E_n}$ is projectively Gorenstein,
because $\omega_{C_{E_n}}$ is trivial, since there is an
everywhere non-zero global section of $\omega_{C_{E_n}}$, given by
the meromorphic 1-form on each component with residues 1 and $-1$
at the nodes (in a suitable order).

All the other $C_{Z_n}$'s, instead, are not locally Gorenstein
because $\omega_{C_{Z_n}}$, although of degree zero, is not
trivial. Indeed a graph $Z_n$, different from $E_n$, certainly has
a vertex with valence 1. This corresponds to a line $l$ such that
$\omega_{C_{Z_n}}\otimes\Oc_l$ is not trivial.
\end{example}

\vspace{-.7cm}

\bigskip

\section{Zappatic surfaces and associated graphs}\label{S:3}

First of all, we need to introduce the singularities we will allow.

\begin{definition}[Zappatic singularity] \label{def:zappaticsing}
Let $X$ be a surface and let $x\in X$ be a point. We will say that
$x$ is a {\it Zappatic singularity} for $X$ if $(X,x)$ is locally
analytically isomorphic to a pair $(Y,y)$ where $Y$ is the cone
over either a curve $C_{T_n}$ or a curve $C_{Z_n}, n\geq 3$, and
$y$ is the vertex of the cone. Accordingly we will say that $x$ is
either a $T_n$- or a $Z_n$-{\em point} for $X$.
\end{definition}

\begin{definition}[Zappatic surface] \label{def:Zappsurf}
Let $X$ be a projective surface with its irreducible components
$X_1,\ldots,X_v$. We will assume that $X$ has the following properties:
\begin{itemize}
\item $X$ is reduced and connected in codimension one;

\item $X_1,\ldots, \; X_v$ are smooth;

\item the singularities in codimension one of $X$ are at most
double curves which are smooth and irreducible;

\item the further singularities of $X$ are Zappatic singularities.
\end {itemize}
A surface like $X$ will be called a {\it Zappatic surface}.
If moreover $X$ is embedded in a projective space $\P ^r$ and
all of its irreducible components are planes,
we will say that $X$ is a {\it planar Zappatic surface}.
\end{definition}

\begin{notation}\label{not:Cij}
Let $X$ be a Zappatic surface.
Let us denote by:
\begin{itemize}
\item $X_i:$ an irreducible component of $X$, $i \leq i \leq v$;
\item $C_{ij}:=X_i \cap X_j$, $1 \leq i \neq j \leq v$, if $X_i$ and $X_j$ meet along a curve,
otherwise set $C_{ij}=\emptyset$;
\item $g_{ij}:$ the genus of $C_{ij}$, $1 \leq i \neq j \leq v$;
\item $C:=\Sing(X)=\cup_{i<j} \, C_{ij}$, the union of all the double curves of $X$;
\item $\Sigma_{ijk}: = X_i \cap X_j \cap X_k$, $1 \leq i \neq j \neq k \leq v$,
if  $X_i \cap X_j \cap X_k\ne\emptyset$, otherwise $\Sigma_{ijk} = \emptyset$;
\item $m_{ijk}:$ the cardinality of the set $\Sigma_{ijk}$;
\item $P_{ijk}^h:$ the Zappatic singular point belonging to $\Sigma_{ijk}$, for $h=1,\ldots,m_{ijk}$.
\end{itemize}
Furthermore, if $X\subset\PR$, for some $r$, we denote by
\begin{itemize}
\item $d :$ the degree of $X$;
\item $d_i :$ the degree of $X_i$, $i \leq i \leq v$;
\item $c_{ij} :$ the degree of $C_{ij}$, $1 \leq i \neq j \leq v$;
\item $D :$ a general hyperplane section of $X$;
\item $g :$ the arithmetic genus of $D$;
\item $D_i :$ the (smooth) irreducible component of $D$ lying in $X_i$,
which is a general hyperplane section of $X_i$, $1 \leq i  \leq v$;
\item $g_i :$ the genus of $D_i$, $1 \leq i \leq v$.
\end{itemize}
If moreover $X$ is a planar Zappatic surface, then $d=v$,
each non-empty set $\Sigma_{ijk}$ is a singleton
and $m_{ijk} =1$, for each $i \neq j \neq k$.
\end{notation}

\begin{remark}\label{rem:pg}
A Zappatic surface $X$ is locally Cohen-Macaulay.
Thus the dualizing sheaf $\omega_X$ is well-defined.
If $X$ has only $E_n$-points as Zappatic singularities, then $X$ is locally Gorenstein,
hence $\omega_X$ is an invertible sheaf.
If $X$ has global normal crossings, i.e.\ if $X$ has only $E_3$-points
as Zappatic singularities,
we define the {\it geometric genus} of $X$ as:
\begin{equation}\label{eq:ggsurfaces}
p_g(X):=h^0(X, \omega_X).
\end{equation}
If $X$ is smoothable, namely if $X$ is the central fibre of a degeneration,
we will define its geometric genus later in Definition \ref{def:pg}.
\end{remark}

\begin{figure}[ht]
\vspace{-24mm}
\[
\begin{array}{cc}
\begin{xy}
<50pt,0pt>:
0; {\xypolygon10{~={18}~:{(1.75,0):}~>{}}},
(0,.4)="a"*=0{\bullet};"7"**@{-}?(.6)*{}="b" ,
"a";"8"**@{-}?(.6)*{}="c" ,
"a";"9"**@{-}?(.6)*{}="d" ,
"b";"c"**@{-}?+<-1pt,-1pt>*!U{\st D_1} ,
"c";"d"**@{-}?+<1pt,-1pt>*!U{\st D_2} ,
"b";"d"**@{.}?(.35)*!D{\st D_3} ,
"7";"8"**@{}?*{X_1} ,
"8";"9"**@{}?*{X_2} ,
"7";"9"**@{}?*+!DR{X_3} ,
"7"*!U{\st C_{13}} ,
"8"*!U{\st C_{12}} ,
"9"*!U{\st C_{23}} ,
\end{xy}
\quad & \quad
\begin{xy}
<50pt,0pt>:
0; {\xypolygon10{~:{(1.75,0):}~>{}}},
(0,.4)="a"*=0{\bullet};"7"**@{-}?(.6)*{}="b" ,
"a";"8"**@{-}?(.6)*{}="c" ,
"a";"9"**@{-}?(.6)*{}="d" ,
"a";"10"**@{-}?(.6)*{}="e" ,
"b";"c"**@{-}?*!DL{\st D_1} ,
"c";"d"**@{-}?+<0pt,4pt>*{\st D_2} ,
"d";"e"**@{-}?*!DR{\st D_3} ,
"7";"8"**@{}?*{X_1} ,
"8";"9"**@{}?*{X_2} ,
"9";"10"**@{}?*{X_3} ,
"8"*!U{\st C_{12}} ,
"9"*!U{\st C_{23}} ,
\end{xy}
\\[33mm]
E_3 \text{-point} & R_3 \text{-point}
\\[-16mm]
\begin{xy}
<50pt,0pt>:
0; {\xypolygon11{~={8.1818}~:{(1.75,0):}~>{}}},
(0,.4)="a"*=0{\bullet};"7"**@{-}?(.6)*{}="b" ,
"a";"8"**@{-}?(.6)*{}="c" ,
"a";"9"**@{-}?(.6)*{}="d" ,
"a";"10"**@{-}?(.6)*{}="e" ,
"a";"11"**@{-}?(.6)*{}="f" ,
"b";"c"**@{-}?+<4pt,2pt>*!D{\st D_1} ,
"c";"d"**@{-}?+<1pt,2pt>*!D{\st D_2} ,
"d";"e"**@{-}?+<-1pt,2pt>*!D{\st D_3} ,
"e";"f"**@{-}?+<-4pt,2pt>*!D{\st D_4} ,
"7";"8"**@{}?*{X_1} ,
"8";"9"**@{}?*{X_2} ,
"9";"10"**@{}?*{X_3} ,
"10";"11"**@{}?*{X_4} ,
"8"+<0pt,-2pt>*!U{\st C_{12}} ,
"9"*!U{\st C_{23}} ,
"10"+<2pt,-2pt>*!U{\st C_{34}} ,
\end{xy}
\quad & \quad
\begin{xy}
<50pt,0pt>:
0; {\xypolygon10{~:{(1.75,0):}~>{}}},
(0,.4)="a"*=0{\bullet};"7"**@{-}?(.6)="b" ,
"a";"8"**@{-}?(.6)="c" ,
"a";"9"**@{-}?(.6)="d" ,
"8";"9"**@{}?="m" ,
"a";"10"**@{-}?(.6)="e" ,
"a";"m"**@{-} ?!{"c";"d"}="n" ,
"b";"c"**@{-}?*!DL{\st D_1} ,
"c";"n"**@{-}?+<0pt,4pt>*{\st D_2} ,
"d";"e"**@{-}?*!DR{\st D_3} ,
"7";"8"**@{}?*{X_1} ,
"8";"m"**@{}?+<3pt,10pt>*!D{X_2} ,
"9";"10"**@{}?*{X_3} ,
"8"*!U{\st C_{12}} ,
"9"+<2pt,0pt>*!U{\st C_{23}} ,
"m";"9"**@{}?(.8)-(0,0.4)="p" ,
"a";"p"**@{-} ?(.6)="q" ,
"m";"9"**@{}?(.35)+<0pt,1pt>*!D{\dt X_4} ,
"m"*!U{\st C_{24}} ,
"n";"q"**@{-}?-<0pt,5pt>*{\st D_4} ,
"n";"d"**@{} ?!{"a";"p"}="r" ,
"n";"r"**@{.} ,
"r";"d"**@{-} ,
\end{xy}
\\[4mm]
R_4 \text{-point} & S_4 \text{-point}
\end{array}
\]
\caption{Examples of good Zappatic singularities.}
\label{fig:R3E3R4}
\end{figure}
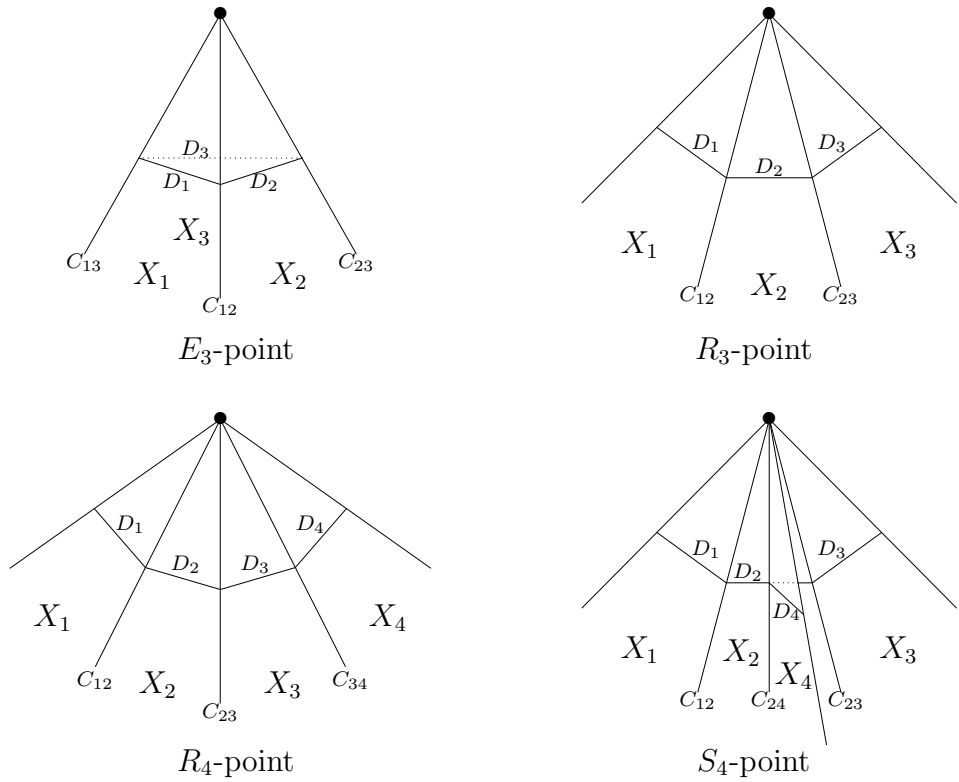

\begin{definition}[Good Zappatic surface] \label{def:goodzapp}
The \emph{good Zappatic singularities} are the
\begin{itemize}
\item $R_n$-points, for $n\ge3$,
\item $S_n$-points, for $n\ge4$,
\item $E_n$-points, for $n\ge3$,
\end{itemize}
which are the Zappatic singularities whose associated stick
curves are respectively $C_{R_n}$, $C_{S_n}$, $C_{E_n}$
(see Examples \ref{ex:tngraphs} and \ref{ex:zngraphs},
Figures \ref{fig:zapsings}, \ref{fig:zapsings2} and \ref{fig:R3E3R4}).

A \emph{good Zappatic surface} is a Zappatic surface
with only good Zappatic singularities.
\end{definition}

To a good Zappatic surface $X$ we can associate a complex $G_X$,
which we briefly call the {\it associated graph} to $X$.

\begin{figure}[ht]
\[
\begin{array}{cccc}
\,\,\,
\begin{xy}
0; <30pt,0pt>: 
(0,1)="a"*=0{\bullet};(1,0)*=0{\bullet}**@{-}="b";(2,1)*=0{\bullet}**@{-}="c";
"a"**@{--} ,
"a"*+!DR{\st v_1} ,
"b"*+!U{\st v_2} ,
"c"*+!DL{\st v_3} ,
\end{xy}
\,\,\, & \,\,\,
\begin{xy}
\xyimport(2,1)(0,0){\rotatebox{45}{\begin{xy} 0; <30pt,0pt>: 
(0,1.4142).(1.4142,0)*[!\colorxy(0.875 0.875 0.875)]\frm{*} , \end{xy}}}
0; <30pt,0pt>: 
(0,1)="a".(2,2)*[white]\frm{*} ,
"a"*=0{\bullet};(1,0)="b"*=0{\bullet}**@{-};(2,1)*=0{\bullet}**@{-}="c";"a"**@{-} ,
"a"*+!DR{\st v_1} , "b"*+!U{\st v_2} , "c"*+!DL{\st v_3} ,
\end{xy}
\,\,\, & \,\,\,
\raisebox{-7.5pt}{$
\begin{xy}
0; <45pt,0pt>: 
(0,1)="a"*=0{\bullet};(0,0)*=0{\bullet}**@{-}="d";(1,0)*=0{\bullet}**@{-}="b";(1,1)*=0{\bullet}**@{-}="c";"a"**@{--} ,
"a"*+!DR{\st v_1} ,
"b"*+!UL{\st v_3} ,
"c"*+!DL{\st v_4} ,
"d"*+!UR{\st v_2} ,
\end{xy}
$}
\,\,\, & \,\,\,
\begin{xy}
0; <30pt,0pt>: 
(0,1)="a"*=0{\bullet};(1,0)*=0{\bullet}**@{-}="b";(2,1)*=0{\bullet}**@{-}="c" ,
"b";(1,1.4142)*=0{\bullet}**@{-}="d" ,
"a"*+!DR{\st v_1} ,
"b"*+!U{\st v_2} ,
"c"*+!DL{\st v_4} ,
"d"*+!D{\st v_3} ,
"b"*\cir<7pt>{ul^dl}*\cir<9pt>{ul^dl} ,
\end{xy}
\,\,\,
\\[2mm]
R_3\text{-point} & E_3\text{-point}
& R_4\text{-point} & S_4\text{-point}
\end{array}
\]
\caption{Associated graphs of $R_3$-, $E_3$-, $R_4$- and $S_4$-points (cf.\ Figure \ref{fig:R3E3R4}).}
\label{fig:graphR3E3R4}
\end{figure}
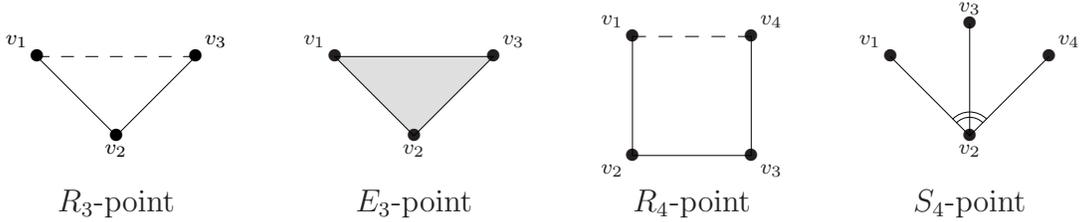

\begin{definition}[The associated graph to $X$] \label{def:dualgraph}
Let $X$ be a good Zappatic surface with Notation \ref{not:Cij}.
The graph $G_X$ associated to $X$ is defined as follows (cf.\ Figure \ref{fig:graphR3E3R4}):

\begin{itemize}

\item each surface $X_i$ corresponds to a vertex $v_i$;

\item each double curve $C_{ij}$ correspond to an edge $e_{ij}$ joining $v_i$ and $v_j$.
The edge $e_{ij}, i<j$, is oriented from the vertex $v_i$ to the one $v_j$;

\item each $E_n$-point $P$ of $X$ is a face of the graph whose $n$
edges correspond to the double curves concurring at $P$. This is
called a {\em $n$-face} of the graph;

\item for each $R_n$-point $P$, with $n\ge 3$, if $P \in
X_{i_1}\cap X_{i_2}\cap\cdots\cap X_{i_n}$, where $X_{i_j}$ meets
$X_{i_k}$ along a curve $C_{i_ji_k}$ only if $1=|j-k|$, we add in
the  graph a \emph{dashed edge} joining the vertices corresponding
to $X_{i_1}$ and $X_{i_n}$. The dashed edge $e_{i_1,i_n}$,
together with the other $n-1$ edges $e_{i_j,i_{j+1}}$,
$j=1,\ldots,n-1$, bound an \emph{open $n$-face} of the graph;

\item for each $S_n$-point $P$, with $n\ge 4$, if $P \in
X_{i_1}\cap X_{i_2}\cap \cdots\cap X_{i_n}$, where $X_{i_1},
\ldots, X_{i_{n-1}}$ all meet $X_{i_n}$ along curves $C_{i_ji_n}$,
$j=1,\ldots,n-1$, concurring at $P$, we mark this in the graph by
an a \emph{angle} spanned by the edges corresponding to the curves
$C_{i_ji_n}$, $j=1,\ldots,n-1$.
\end{itemize}

In the sequel, when we speak of {\it faces} of $G_X$ we always mean closed faces.
Of course each vertex $v_i$ is weighted with the relevant
invariants of the corresponding surface $X_i$.
We will usually omit these weights if $X$ is planar,
i.e.\ if all the $X_i$'s are planes.
\end{definition}

Since each $R_n$-, $S_n$-, $E_n$-point is an element of some set of points $\Sigma_{ijk}$
(cf.\ Notation \ref{not:Cij}),
we remark that there can be different
faces (as well as open faces and angles) of $G_X$ which
are incident on the same set of vertices and edges.
However this cannot occur if $X$ is planar.

\vspace{-5mm}

Notice that angles, open and closed faces of $G_X$ have been defined
in order to encode the good Zappatic singularities of $X$.
In other words, the associated graph $G_X$ uniquely determines
the configuration of the good Zappatic singularities of $X$.

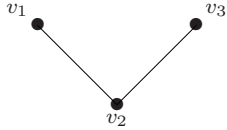
\begin{figure}[ht]
\[
\begin{xy}
0; <30pt,0pt>: 
(0,1)="a"*=0{\bullet};(1,0)*=0{\bullet}**@{-}="b";(2,1)*=0{\bullet}**@{-}="c",
"a"*+!DR{\st v_1} ,
"b"*+!U{\st v_2} ,
"c"*+!DL{\st v_3} ,
\end{xy}
\]
\caption{Associated graph of a $R_3$-point in a good, planar Zappatic surface.} \label{fig:graphR3plan}
\end{figure}

Consider three vertices $v_i, v_j, v_k$ of $G_X$ in such a way
that $v_i$ is joint with $v_j$ and $v_k$. Any point in
$C_{ij}\cap C_{ik}$ is either a $R_n$-, or a $S_n$-, or an
$E_n$-point, and the curves $C_{ij}$ and $C_{ik}$ intersect
transversally, by definition of Zappatic singularities.
Hence we can compute the intersection number $C_{ij}\cdot C_{ik}$
by adding the number of closed and open faces and of angles
involving the edges $e_{ij}, e_{ik}$. In particular, if $X$ is
planar, for every pair of adjacent edges only one of
the following possibilities occur: either they belong to an open
face, or to a closed one, or to an angle. Therefore for good, planar Zappatic surfaces
we can avoid to mark open $3$-faces without loosing any information
(see Figure \ref{fig:graphR3plan}, cf.\ Figure \ref{fig:graphR3E3R4}).

\vspace{-2mm}

\begin{remark}\label{rem:orientation}
We also notice that, given our choices, if $X$ is good Zappatic
and has only $E_3$-points, the graph $G_X$ comes with a
lexicographic orientation of the faces; indeed, let
$\Sigma_{ijk}=X_i\cap X_j\cap X_k = \{P_{ijk}^1, P_{ijk}^2,
\ldots, P_{ijk}^{m_{ijk}} \}$; thus, each face of $G_X$
corresponds to a sequence of three vertices $i,j,k$ with $i<j<k$,
together with an integer $t$ such that $1 \leq t \leq m_{ijk}$,
hence it will be denoted by $f^t_{\sigma(i),\sigma(j),\sigma(k)}$,
with $\sigma$ any permutation of $i,j,k$, and will be oriented
according to the orientation of its boundary determined by the
sequence of vertices $v_i,v_j,v_k$.

\end{remark}

As for stick curves, if $G$ is a given graph as above,
there does not necessarily exist a good planar
Zappatic surface $X$ such that its associated graph is $G=G_X$.

\begin{figure}[h]
\[
\begin{xy}
0; <45pt,0pt>: 
(0,1)="a", (1,0)="b" ,
"a"."b"*[!\colorxy(0.875 0.875 0.875)]\frm{**} ,
"a"*=0{\bullet};(0,0)*=0{\bullet}**@{-}="d";
"b"*=0{\bullet}**@{-};(1,1)*=0{\bullet}**@{-}="c";
"a"**@{-};"b"**@{-} ,
"a"*+!DR{\st v_1} ,
"b"*+!UL{\st v_3} ,
"c"*+!DL{\st v_4} ,
"d"*+!UR{\st v_2} ,
\end{xy}
\]
\caption{Graph associated to an impossible planar Zappatic surface.} \label{fig:impo}
\end{figure}
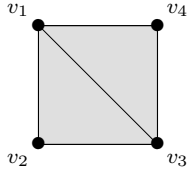

\begin{example}\label{ex:impo}
Consider the graph $G$ of Figure \ref{fig:impo}.
If $G$ were the associated graph of a good planar Zappatic surface
$X$, then $X$ should be a global normal crossing union
of $4$ planes with $5$ double lines and two $E_3$ points,
$P_{123}$ and $P_{134}$, both lying on the double line $C_{13}$.
Since the lines $C_{23}$ and $C_{34}$ (resp.\ $C_{14}$ and
$C_{12}$) both lie on the plane $X_3$ (resp.\ $X_1$), they should
intersect. This means that the planes $X_2,X_4$ also should
intersect along a line, therefore the edge $e_{24}$ should appear in the graph.
\end{example}

\vspace{-3mm}

\medskip

Before going on, we need some notation.

\begin{notation}\label{def:vefg}
Let $X$ be a good Zappatic surface (with invariants as in
Notation \ref{not:Cij}) and let $G= G_X$ be its associated graph. We denote by
\begin{itemize}

\item $v :$ the number of vertices of $G$ (i.e.\ the number
of irreducible components of $X$);

\item $V :$ the (indexed) set of vertices of $G$;

\item $e :$ the number of edges of $G$
(i.e.\ the number of double curves in $X$);

\item $E :$ the set of edges of $G$; this is indexed by the ordered pairs $(i,j) \in V \times V$,
$i <j$, such that the corresponding surfaces $X_i$, $X_j$ meet along the
curve $C_{ij} = C_{ji}$;

\item $f_n :$ the number of $n$-faces of
$G$, i.e.\ the number of $E_n$-points of $X$, for $n \geq 3$;

\item $f := \sum_{n \geq 3} f_n$, the number of faces of $G$,
i.e.\ the total number of $E_n$-points of $X$, for all $n \geq3$;

\item $r_n :$ the number of open $n$-faces of $G$, i.e.\ the number of
$R_n$-points of $X$, for $n\ge 3$;

\item $s_n :$ the number of $n$-angles of $G$, i.e.\ the number of
$S_n$-points of $X$, for $n\ge 4$;

\item $\chi(G) :=v-e+f$, i.e.\ the Euler-Poincar\'e characteristic of $G$;

\item $G^{(1)} :$ the \emph{$1$-skeleton} of $G$, i.e.\ the graph
obtained from $G$ by forgetting all the faces, dashed edges and angles;

\item $\chi(G^{(1)}) =v-e$, i.e.\ the Euler-Poincar\'e characteristic of $G^{(1)}$.
\end{itemize}
\end{notation}

\begin{remark}\label{rem:G^(1)}
Observe that, when $X$ is a good, planar Zappatic surface, the $1$-skeleton $G^{(1)}_X$ of $G_X$
coincides with the dual graph $G_D$ of the general hyperplane section $D$ of $X$.
\end{remark}

Now we can compute some of the invariants of good Zappatic surfaces.

\begin{proposition}\label{prop:ghypsect}
Let $X = \bigcup_{i=1}^v X_i \subset \PR$ be a good Zappatic
surface and let $G = G_X$ be its associated graph. Let $C$ be the
double locus of $X$, i.e.\ the union of the double curves of $X$,
$C_{ij} = C_{ji} = X_i \cap X_j$ and let $c_{ij} = \deg(C_{ij})$.
Let $D_i$ be a general hyperplane section of $X_i$, and denote by $g_i$ its genus. Then
the arithmetic genus of a general hyperplane section $D$ of $X$ is:
\begin{equation}
g = \sum_{i=1}^v g_i + \sum_{e_{ij}\in E} \; c_{ij} -v +1.
\label{eq:g}
\end{equation}
In particular, when $X$ is a good, planar Zappatic surface, then
\begin{equation}\label{eq:gplanar}
g = e- v + 1 = 1 - \chi(G^{(1)}).
\end{equation}
\end{proposition}
\begin{proof} Denote by $d_i$ the degree of $X_i$, $1 \leq i \leq v$. Then,
$D$ is the union of the $v$ irreducible components $D_i$, $1 \leq
i \leq v$, such that $\deg(D_i) = d_i$ and $d := \deg(D) =
\sum_{i=1}^v d_i$. Consider its associated graph $G_D$, defined as
in \S \ref{S:2}.

Take $G$, whose indexed set of edges is denoted by $E$, and
consider $e_{ij} \in E$ joining its vertices $v_i$ and $v_j$, $i
<j$, which correspond to the irreducible components $X_i$ and
$X_j$, respectively. Since $e_{ij}$ in $G$ correspond to the
double curve $C_{ij}$, we have exactly $c_{ij}$ oriented edges in
the graph $G_D$ joining its vertices $v_i$ and $v_j$, which now
correspond to the irreducible components of $D$, $D_i$ and $D_j$,
respectively. These $c_{ij}$ oriented edges correspond to the
$c_{ij}$ nodes of the reducible curve $D_i \cup D_j$, which is
part of the hyperplane section $D$.

Now, recall that the Hilbert polynomial of $D$ is, with our
notation, $P_{D}(t) = d t + 1 - g$. On the other hand, $P_{D}(t)$
equals the number of independent conditions imposed on
hypersurfaces $\Ha$ of degree $t \gg 0$ to contain $D$.

From what observed above on $G_D$, it follows that the number of singular
points of $D$ is $\sum_{e_{ij}\in E} c_{ij}$. These points impose independent
conditions on hypersurfaces $\Ha$ of degree $t \gg 0$.

Since $t\gg0$ by assumption, we get that the map
$$
H^0(\Oc_{\PR} (t)) \to H^0(\Oc_{D_i}(t))
$$
is surjective and that the line bundle $\Oc_{D_i}(t)$ is
non-special on $D_i$, for each $1 \leq i \leq v$. Thus, in order
for $\Ha$ to contain $D_i$ we have to impose $d_i t - g_i + 1 -
\sum_{j \; {\rm s.t.} \;e_{ij} \in E } c_{ij}$ conditions. Therefore the total
number of conditions for $\Ha$ to contain $D$ is:
\begin{align*}
\sum_{e_{ij}\in E}  c_{ij} + \sum_{i=1}^v \biggl( d_i t -g_i +1 -
\sum_{j, e_{ij} \in E} c_{ij} \biggr) &= \sum_{e_{ij}\in E} c_{ij}
+ dt - \sum_{i=1}^v g_i + v
 - \sum_{i=1}^v \sum_{j, e_{ij} \in E } c_{ij} = \\
&= d t + v - \sum_{i=1}^v g_i - \sum_{e_{ij}\in E} c_{ij},
\end{align*}
since $\sum_{i=1}^v \sum_{j, e_{ij} \in E } \; c_{ij} = 2 \sum_{e_{ij}\in E} \; c_{ij}$.
This proves \eqref{eq:g} (cf.\ formula \eqref{eq:genuscurves}).

The second part of the statement directly follows from the above
computations and from the fact that, in the good planar Zappatic
case $g_i=0$ and $c_{ij} = 1$, for each $i < j$, i.e.\ $G_D$ coincides with
$G^{(1)}$ (cf.\ Remark \ref{rem:G^(1)}).
\end{proof}

By recalling Notation \ref{def:vefg}, one also has:

\begin{proposition}\label{prop:4.chi}
Let $X = \bigcup_{i=1}^v X_i $ be a good Zappatic surface and
$G_X$ be its associated graph. Let $C$ be the double locus of $X$, which is the
union of the curves $C_{ij} = C_{ji} = X_i \cap X_j$. Then:
\begin{equation} \label{eq:chi}
\chi(\Oc_X) = \sum_{i=1}^v \chi(\Oc_{X_i}) - \sum_{e_{ij} \in E}
\chi(\Oc_{C_{ij}}) + f .
\end{equation}

In particular, when $X$ is a good, planar Zappatic surface, then
\begin{equation} \label{eq:chiplan}
\chi(\Oc_X) = \chi(G_X) = v -e + f .
\end{equation}
\end{proposition}
\begin{proof}

We can consider the sheaf morphism:
\begin{equation}\label{eq:aiuto2}
\bigoplus_{i=1}^v \Oc_{X_i}
\xrightarrow{\;\lambda\;} \bigoplus_{1 \le i<j \le v} \Oc_{C_{ij}} ,
\end{equation}
defined in the following way: if
$$
\pi_{ij}: \bigoplus_{1 \le i < j \le v} \Oc_{C_{ij}} \to
\Oc_{C_{ij}}
$$
denotes the projection on the $(ij)^{\rm th}$-summand, then
$$
(\pi_{ij} \circ \lambda) (h_1, \ldots , h_v) := h_i - h_j.
$$
Notice that the definition of $\lambda$ is consistent with the
lexicographic order of the indices and with the lexicographic
orientation of the edges of the graph $G_X$.

Observe that, if $\tilde{X}$ denotes the desingularization of $X$,
then $\tilde{X}$ is isomorphic to the disjoint union of the
smooth, irreducible components $X_i$, $1 \leq i \leq v$, of $X$.
Therefore, by the very definition of $\Oc_X$, we see that
$$
\ker(\lambda) \cong \Oc_X.
$$

We claim that the morphism $\lambda$ is not surjective and that
its cokernel is a sky-scraper sheaf supported at the $E_n$-points
of $X$. To show this, we focus on any irreducible component of $C
= \bigcup_{1 \leq i < j \leq v} C_{ij}$, the double locus of $X$.

Fix any index pair $(i,j)$, with $i <j$, and consider the
generator
\begin{equation}\label{eq:gener}
(0, \ldots, 0 , 1 , 0 , \ldots, 0) \in \bigoplus_{1 \leq l < m
\leq v} \Oc_{C_{lm}},
\end{equation}
where $1 \in \Oc_{C_{ij}}$, the $(ij)^{\rm th}$-summand. The
obstructions to lift up this element to an element of
$\bigoplus_{1 \leq t \leq v} \Oc_{X_t}$ are given by the presence
of good Zappatic singularities of $X$ along $C_{ij}$.

For what concerns the irreducible components
of $X$ which are not involved in the intersection determining a
good Zappatic singularity on $C_{ij}$, the element in
\eqref{eq:gener} trivially lifts-up to $0$ on each of them.
Thus, in the sequel, we shall focus only on the irreducible components
involved in the Zappatic singularity, which will be denoted by $X_i, \; X_j, \;
X_{l_t}$, for $ 1 \leq t \leq n-2$.

We have to consider different cases, according to the good
Zappatic singularity type lying on the curve $C_{ij} = X_i \cap
X_j$.

\begin{itemize}
\item Suppose that $C_{ij}$ passes through a $R_n$-point $P$ of
$X$, for some $n$; we have two different possibilities. Indeed:

\noindent (a) let $X_i$ be an ``external" surface for $P$ ---
i.e.\ $X_i$ corresponds to a vertex of the associated graph of $P$
which has valence $1$. Therefore, we have:
\[
\begin{picture}(150,15)
\thinlines \put(0,10){\line(1,0){75}}
\put(150,10){\line(-1,0){45}} \put(0,10){\circle*{3}}
\put(30,10){\circle*{3}} \put(60,10){\circle*{3}}
\put(120,10){\circle*{3}} \put(150,10){\circle*{3}}
\put(84,7.7){$\cdots$} \put(-5,2){$\scriptstyle X_i$}
\put(25,2){$\scriptstyle X_j$} \put(55,2){$\scriptstyle X_{l_1}$}
\put(110,2){$\scriptstyle X_{l_{n-3}}$} \put(140,2){$\scriptstyle
X_{l_{n-2}}$}
\end{picture}
\]
%
In this situation, the element in \eqref{eq:gener} lifts up to
$$
(1, 0, \ldots, 0) \in \Oc_{X_i} \oplus \Oc_{X_j} \oplus
\bigoplus_{1\leq t \leq n-2} \Oc_{X_{l_t}}.
$$

\noindent (b) let $X_i$ be an ``internal" surface for $P$ ---
i.e.\ $X_i$ corresponds to a vertex of the associated graph to $P$
which has valence $2$. Thus, we have a picture like:
\[
\begin{picture}(210,15)
\thinlines \put(0,10){\line(1,0){135}}
\put(210,10){\line(-1,0){45}} \put(0,10){\circle*{3}}
\put(30,10){\circle*{3}} \put(60,10){\circle*{3}}
\put(90,10){\circle*{3}} \put(120,10){\circle*{3}}
\put(180,10){\circle*{3}} \put(210,10){\circle*{3}}
\put(144,7.7){$\cdots$} \put(-5,2){$\scriptstyle {X_{l_1}}$}
\put(25,2){$\scriptstyle {X_{l_2}}$} \put(55,2){$\scriptstyle
{X_{l_3}}$} \put(85,2){$\scriptstyle {X_i}$}
\put(115,2){$\scriptstyle {X_j}$} \put(170,2){$\scriptstyle
{X_{l_{n-3}}}$} \put(200,2){$\scriptstyle {X_{l_{n-2}}}$}
\end{picture}
\]

In this case, the element in \eqref{eq:gener} lifts up to the
$n$-tuple having components:
\begin{itemize}
\item[] $1 \in \Oc_{X_i}$, \item[] $0 \in \Oc_{X_j}$, \item[] $1
\in \Oc_{X_{l_t}}$, for those $X_{l_t}$'s corresponding to
vertices in the graph associated to $P$ which are on the left of
$X_i$ and, \item[] $0 \in \Oc_{X_{l_k}}$ for those
$X_{l_k}$'s corresponding to vertices in the graph
associated to $P$ which are on the right of $X_j$.
\end{itemize}
\item Suppose that $C_{ij}$ passes through a $S_n$-point $P$ of
$X$, for any $n$; as before, we have two different possibilities.
Indeed:

\noindent (a) let $X_i$ corresponds to the vertex of valence $n-1$
in the associated graph to $P$, i.e.
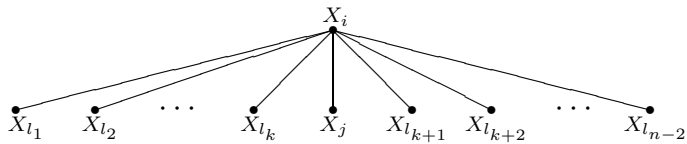
\[
\begin{picture}(240,50)
\thinlines \put(120,40){\line(0,-1){30}}
\put(120,40){\line(-1,-1){30}} \put(120,40){\line(-3,-1){90}}
\put(120,40){\line(-4,-1){120}} \put(120,40){\line(1,-1){30}}
\put(120,40){\line(2,-1){60}} \put(120,40){\line(4,-1){120}}
\put(120,40){\circle*{3}} \put(0,10){\circle*{3}}
\put(30,10){\circle*{3}} \put(90,10){\circle*{3}}
\put(120,10){\circle*{3}} \put(150,10){\circle*{3}}
\put(180,10){\circle*{3}} \put(240,10){\circle*{3}}
\put(204,7.7){$\cdots$} \put(54,7.7){$\cdots$}
\put(116,43){$\scriptstyle {X_i}$} \put(-3,2){$\scriptstyle
{X_{l_1}}$} \put(25,2){$\scriptstyle {X_{l_2}}$}
\put(85,2){$\scriptstyle {X_{l_k}}$} \put(115,2){$\scriptstyle
{X_j}$} \put(140,2){$\scriptstyle {X_{l_{k+1}}}$}
\put(170,2){$\scriptstyle {X_{l_{k+2}}}$}
\put(230,2){$\scriptstyle {X_{l_{n-2}}}$}
\end{picture}
\]

In this situation, the element in \eqref{eq:gener} lifts up to the
$n$-tuple having components:
\begin{itemize}
\item[] $1 \in \Oc_{X_i}$, \item[] $0 \in \Oc_{X_j}$, \item[] $1
\in \Oc_{X_{l_t}}$, for all $1 \leq t \leq n-2$.
\end{itemize}

\noindent (b) let $X_i$ corresponds to a vertex of valence $1$
in the associated graph to $P$. Since $C_{ij} \neq \emptyset$ by
assumption, then $X_j$ has to be the vertex of valence $n-1$, i.e.
we have the following picture:
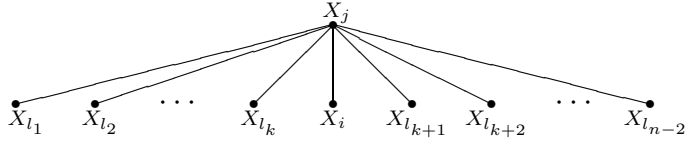
\[
\begin{picture}(240,50)
\thinlines \put(120,40){\line(0,-1){30}}
\put(120,40){\line(-1,-1){30}} \put(120,40){\line(-3,-1){90}}
\put(120,40){\line(-4,-1){120}} \put(120,40){\line(1,-1){30}}
\put(120,40){\line(2,-1){60}} \put(120,40){\line(4,-1){120}}
\put(120,40){\circle*{3}} \put(0,10){\circle*{3}}
\put(30,10){\circle*{3}} \put(90,10){\circle*{3}}
\put(120,10){\circle*{3}} \put(150,10){\circle*{3}}
\put(180,10){\circle*{3}} \put(240,10){\circle*{3}}
\put(204,7.7){$\cdots$} \put(54,7.7){$\cdots$}
\put(116,43){$\scriptstyle {X_j}$} \put(-3,2){$\scriptstyle
{X_{l_1}}$} \put(25,2){$\scriptstyle {X_{l_2}}$}
\put(85,2){$\scriptstyle {X_{l_k}}$} \put(115,2){$\scriptstyle
{X_i}$} \put(140,2){$\scriptstyle {X_{l_{k+1}}}$}
\put(170,2){$\scriptstyle {X_{l_{k+2}}}$}
\put(230,2){$\scriptstyle {X_{l_{n-2}}}$}
\end{picture}
\]
Thus, the element in \eqref{eq:gener} lifts up to the $n$-tuple
having components
\begin{itemize}
\item[] $1 \in \Oc_{X_i}$, \item[] $0 \in \Oc_{X_j}$, \item[] $0
\in \Oc_{X_{l_t}}$, for all $1 \leq t \leq n-2$.
\end{itemize}

\item Suppose that $C_{ij}$ passes through an $E_n$-point $P$ for
$X$. Then, each vertex of the associated graph to $P$ has valence
$2$. Since such a graph is a cycle, it is clear that no lifting of
\eqref{eq:gener} can be done.
\end{itemize}

To sum up, we see that $\coker (\lambda)$ is supported at the
$E_n$-points of $X$. Furthermore, if we consider
\[
\bigoplus_{1 \le i < j \le v} \Oc_{C_{ij}}
\xrightarrow{\;ev_P\;} \Oc_P = \underline{\CC}_P, \qquad \oplus
f_{ij} \mapsto \sum f_{ij} (P)
\]
it is clear that, if $P$ is an $E_n$-point then
$$
ev_P \left(\bigoplus_{1 \le i < j \le v} \Oc_{C_{ij}} /
\im(\lambda)\right) \cong \underline{\CC}_P.
$$
This means that
$$\coker (\lambda) \cong \underline{\CC}^f.$$

By the exact sequences
\[
0 \to \Oc_X  \to  \bigoplus_{1 \le i \le v} \Oc_{X_i}  \to
\im(\lambda)  \to 0, \qquad 0 \to \im (\lambda) \to  \bigoplus_{1
\le i < j \le v} \Oc_{C_{ij}}  \to \underline{\CC}^f  \to 0,
\]
we get \eqref{eq:chi}.
\end{proof}

In the next section we will see how the computation
of the geometric genus is much more involved, even in the case of $X$ having
only $E_3$-points as Zappatic singularities.

\bigskip

\section{The geometric genus of a Zappatic surface with only $E_3$-points}\label{S:4pg}

The main purpose of this section is to
compute the {\em geometric genus} of a projective good Zappatic surface
$X = \bigcup_{i=1}^v X_i$ (cf.\ Remark \ref{rem:pg}), which is assumed to have only
$E_3$-points, i.e.\ global normal crossing singularities.
In terms of its associated graph, this means that $G_X$ is a
subgraph of the complete graph on $v$ vertices, which has only
$3$-faces (i.e.\ triangles).

We first want to recall some definitions and results, which will
be used in the sequel.

\begin{definition}\label{def:poincare}
Let $T$ be a smooth surface and $C$ be a smooth, irreducible curve
on $T$. Let $\omega$ be a global meromorphic 2-form on $T$ whose
polar locus contains $C$. We may assume that $x,y$ are local
coordinates on $T$ in an analytic neighbourhood of a point of $C$
in such a way that $y=0$ is the local equation defining $C$. Then,
the \emph{Poincar\'e residue map} (or {\em adjunction map})
$$\omega_T\otimes\Oc_T(C)\xrightarrow{\,\R_C\,} \omega_C$$
is locally defined by:
\[
\omega=\frac{f(x,y)}{y} \;dx\wedge dy\ \mapsto\  - f(x,0) \;dx
\]
and $- f(x,0) \;dx$ is called the \emph{(Poincar\'e) residue} of
the 2-form $\omega$ along $C$ (see \cite[p.\ 147]{GH}).

If, more generally, $C$ is assumed to be a reduced (possibly
reducible) divisor with only normal crossing singularities, denote
by $\omega_C$ its dualizing sheaf and take local coordinates on
$T$ in an analytic neighbourhood of a node of $C$ in such a way
that $xy=0$ is its local defining equation. If
$$\omega = \frac{f(x,y)}{xy} \;dx\wedge dy,$$
then one defines the Poincar\'e residue map by considering
\begin{equation}\label{eq:presnod}
\R_C(\omega) \in H^0(C, \omega_C)
\end{equation}defined as the pair of forms:

\begin{itemize}
\item[(i)] $\displaystyle\omega_x = - \frac{f(x,0)}{x} \;dx$ on
the branch $y=0$, \item[(ii)] $\displaystyle\omega_y=
\frac{f(0,y)}{y} \;dy$ on the branch $x =0$.
\end{itemize}
\end{definition}

\begin{remark}\label{rem:residuinodi}
If $\R_0(\omega_x)$ denotes the usual Poincar\'e residue at the
point $x=0$ of the meromorphic 1-form $\omega_x$ on the (smooth)
branch $y=0$, observe that
\begin{equation}\label{eq:residuinodi}
\R_0(\omega_y) = -\R_0 (\omega_x).
\end{equation}
This is consistent with the definition
of $H^0( C, \omega_C)$. Indeed, assume that $C$ has only $m$ nodes;
then, if $\nu: \tilde{C} \to C$ is the normalization morphism and
if $\{ q_1, q_1' \}, \; \{ q_2, q_2' \}, \ldots, \; \{ q_m, q_m'
\}$ are the pre-images in $\tilde{C}$ of the $m$ nodes of $C$,
then $\omega_C$ is the invertible subsheaf
$$ \omega_C \subset \nu_* (\omega_{\tilde{C}} (\sum_{i=1}^m (q_i + q_i')))$$
such that a section $ \sigma$ of $\nu_* (\omega_{\tilde{C}}
(\sum_{i=1}^m (q_i + q_i')))$, viewed as a section of
$\omega_{\tilde{C}} (\sum_{i=1}^m (q_i + q_i'))$, is a section of
$\omega_C$ if and only if $$\R_{q_i} (\sigma) + \R_{q_i'} (\sigma)
= 0 , \; \; 1 \leq i \leq m.$$
\end{remark}

Unless otherwise stated, from now on $X=\bigcup_{i=1}^v X_i$ will
denote a projective, Zappatic surface with $E_3$-points only as
Zappatic singularities and we use notation as in Definition
\ref{def:Zappsurf} and in Notation \ref{def:vefg}.

It is well known that, for each $i$:
\begin{equation}\label{eq:K_Si}
\omega_{X_i} \cong \omega_X \otimes \Oc_{X_i}(-C_i),\qquad
\text{with } C_i := X_i \cap \overline{(X \setminus X_i)} =
\sum_{j\ne i} C_{ij},
\end{equation}
where $\omega_{X_i}$ is the canonical line bundle of $X_i$,
whereas $\omega_X$ denotes the dualizing sheaf of $X$. Note that
$\omega_X$ is an invertible sheaf by the hypotheses on $X$.

As in Remark \ref{rem:pg}, recall that
the geometric genus of $X$ is denoted
by $p_g(X)$ and defined as $p_g(X) = h^0 (X, \omega_X)$. In order
to compute $p_g(X)$, we need some further remarks which will be
fundamental in the sequel.

\begin{remark}\label{rem:indices}
Observe that, if $X=\bigcup_{i=1}^v X_i$ is as above, since the
intersection $X_i \cap X_j$ --- when non-empty --- is the double curve
$C_{ij} = C_{ji}$, the index pair $(i,j)$ with $i<j$ uniquely
determines the double curve $C_{ij}$. In the same way, when $ i <
j$, since the intersection $X_i \cap X_j \cap X_k$
--- when non-empty ---
is the triple point set $\Sigma_{ijk} = \Sigma_{\tau(i)
\tau(j) \tau(k)}$, for $ k \neq i,j$ and for any $\tau \in \Sym(\{
i, j, k \})$, then the lexicographically ordered index triple
uniquely determines the triple point set either $\Sigma_{kij}$, or
$\Sigma_{ikj}$, or $\Sigma_{ijk}$, according to the fact that
either $k < i$, or $ i < k <j$, or $k > j$.
\end{remark}

\begin{remark}\label{rem:planes}
Since the only Zappatic singularities of $X$ are assumed to be
$E_3$-points, then $G_X$ contains neither dashed edges, nor
angles, nor open faces, nor $n$-faces, with $n\geq 4$. Furthermore,
in such a case the graph $G_X$ comes with a lexicographic orientation
of the faces (see Remark \ref{rem:orientation}). If $X$ is, in particular, 
planar recall that we have strong constraints
on the possible shape of the graph $G_X$ --- because of the geometry of
planes (cf.\ Example \ref{ex:impo})--- and each non-zero $m_{ijk}$ equals one.
\end{remark}

Observe that by the connectedness hypothesis
of $X=\bigcup_{i=1}^v X_i$, we get that $C_i \neq \emptyset$, for
each $ 1 \leq i \leq v$. For simplicity of notation, in the sequel
we shall always denote by $C_{ij}$ the intersection of $X_i$ and
$X_j$, for any $ 1 \leq i < j \leq v$, with the obvious further
condition that $C_{ij}= C_{ji} = \emptyset$ when the index pair
corresponds to two disjoint surfaces in $X$, i.e.\ when there is no
edge $(v_i, v_j)$ in the associated graph $G_X$.

We can define a natural map
\begin{equation}\label{eq:4.fi}
\Phi: \bigoplus_{i=1}^v H^1(X_i, \Oc_{X_i})\to \bigoplus_{1 \le i
< j \le v} H^1(C_{ij}, \Oc_{C_{ij}})
\end{equation}in the following way: if
$$\pi_{ij}: \bigoplus_{1 \le i < j \le
v} H^1(C_{ij}, \Oc_{C_{ij}}) \to H^1(C_{ij},
\Oc_{C_{ij}})$$denotes the projection on the $(ij)^{\rm th}$-summand
and if $r^{(i)}_{C_{ij}}: H^1(X_i, \Oc_{X_i}) \to H^1 (C_{ij},
\Oc_{C_{ij}}) $ denotes the natural restriction map to $C_{ij}$ as
a divisor in $X_i$, where $i <j$, then
\begin{equation}\label{eq:fi}
(\pi_{ij} \circ \Phi) ((a_1, \ldots, a_v)) :=
r^{(i)}_{C_{ij}}(a_i) - r^{(j)}_{C_{ij}}(a_j).
\end{equation}

\begin{remark}\label{rem:segnocurve}
Observe that the above definition is consistent with the
lexicographic order of the indices $1 \leq i \leq v$. In other words,
\eqref{eq:fi} means that we consider the curve $C_{ij}$ as a {\em
positive curve} on the surface $X_i$ and as a {\em negative curve}
on the surface $X_j$, when $i < j$. Furthermore, when the index
pair is such that $C_{ij} = \emptyset$, we obviously consider
$H^1(C_{ij}, \Oc_{C_{ij}})$ as the zero-vector space and $\pi_{ij}
\circ \Phi$ as the zero-map.
\end{remark}

Take an index pair $i <j$ such that $C_{ij} \neq \emptyset$. By
the adjunction sequence of $C_{ij}$ on $X_i$ and on $X_j$, we can
consider the two obvious coboundary maps:
\begin{equation}\label{eq:fid}
\xymatrix@R-5mm@C+3mm{%
& H^1(X_i, \omega_{X_i}) \\
H^0(C_{ij}, \omega_{C_{ij}}) \ar[]!UR;[ur]!DL^{\delta_i} \ar[]!DR;[dr]!UL^{\delta_j} \\
 & H^1(X_j, \omega_{X_j}).
}
\end{equation}
On the other hand, when the index pair $i < j$ is such that
$C_{ij} = \emptyset$, then $H^0(C_{ij}, \omega_{C_{ij}})$ is
considered as the zero-vector space and \eqref{eq:fid} are
the zero-maps. Then, we can define the map
\begin{equation}\label{eq:4.fid}
\Delta: \bigoplus_{1 \le i < j \le v} H^0(C_{ij}, \omega_{C_{ij}})
\to \bigoplus_{i=1}^v H^1(X_i, \omega_{X_i})
\end{equation}in the following way: if
$$\iota_{ij}: H^0(C_{ij}, \omega_{C_{ij}}) \hookrightarrow \bigoplus_{1 \le i <
j \le v} H^0(C_{ij}, \omega_{C_{ij}})$$denotes the natural
inclusion of the $(ij)^{\rm th}$-summand and if $\gamma_{ij} \in
H^0(C_{ij}, \omega_{C_{ij}})$, then
\begin{equation}\label{eq:fid2}
(\Delta \circ \iota_{ij})(\gamma_{ij}) := (0, \ldots, 0,
\delta_i(\gamma_{ij}), 0, \ldots, 0, - \delta_j(\gamma_{ij}), 0 ,
\ldots, 0),
\end{equation}where $i<j$.

Observe that the definition of $\Delta$ is consistent with the
lexicographic order of the indices $1 \leq i \leq v$ and with our
Remark \ref{rem:segnocurve}.

The following preliminary result is an obvious
consequence of our definitions.

\begin{proposition}\label{prop:Fi}
With notation as above, we have
$$\Delta = \Phi^{\vee}.$$
\end{proposition}
\begin{proof}
The proof directly follows from Serre's duality on each summand
and from the fact that the matrix which represents $\Delta$ is the
transpose of the one representing $\Phi$.
\end{proof}

We are now able to prove the main result of this section.

\begin{theorem}\label{thm:4.pgbis}
Let $X=\bigcup_{i=1}^v X_i$ be a projective, good Zappatic surface
with only $E_3$-points as Zappatic singularities. Denote by
$\omega_X$ the dualizing sheaf of $X$ and by $G_{X}$ the
associated graph of $X$ (see Definition \ref{def:dualgraph}). Let
$\Phi$ be the map defined in \eqref{eq:4.fi} and let $p_g(X)$ be
the geometric genus of $X$ as in Remark \ref{rem:pg}. Then the
following inequality:
\begin{equation}\label{eq:Fpgbis}
p_g(X) \leq b_2(G_X)+\sum_{i=1}^v p_g(X_i)+ \dim(\coker(\Phi))
\end{equation}holds, where as costumary $b_2(G_X)$ is the second Betti number of
$G_X$.

Furthermore, a sufficient condition for the equality in
\eqref{eq:Fpgbis} to hold is either that

\begin{itemize}

\item[(i)] each irreducible component $X_i$ is a regular surface,
for $1 \leq i \leq v$, or that
\item[(ii)] for any irregular component $X_j$ of $X$,
the divisor $C_{j} = X_{j} \cap \overline{(X \setminus X_{j})}$ is
ample on $X_j$.

\end{itemize}
\end{theorem}
\begin{proof}To prove the first part of the statement,
we construct a homomorphism
\begin{equation}\label{eq:kerfbis}
H^0(\omega_X) \mapright{f} H_2(G_X, \CC)
\end{equation}and we show that
\begin{equation}\label{eq:kerf2bis}
\ker(f) \cong \bigoplus_{i=1}^v H^0 (X_i , \omega_{X_i}) \oplus
(\coker (\Phi)).
\end{equation}
Then, for what concerns the second part, we prove that either
hypothesis $(i)$ or hypothesis $(ii)$ implies the surjectivity of
the map \eqref{eq:kerfbis}.

Recall that, from our hypotheses it follows that $G_X$ is a
subgraph of the complete graph on $d$ vertices which contains only
$3$-faces and that there can be more than one $3$-face incident
on the same triple of vertices (equiv. edges).
This occurs when (cf.\ Notation \ref{not:Cij}) $m_{ijk} > 1$, for a given triple $v_i,
v_j, v_k$ of vertices of $G_X$. It is obvious that, when
$\Sigma_{ijk}  = \emptyset$, then $m_{ijk} = 0$.

To construct $f$, by \eqref{eq:K_Si}, we consider a global section
$\omega \in H^0(X, \omega_X)$ as a collection
$$\{ \omega_i\}_{1 \leq i \leq v} \in \bigoplus_{1 \leq i \leq v}
H^0(X_i, \omega_{X_i}(C_i)),$$where each $\omega_i$ is a global
meromorphic $2$-form on the corresponding irreducible component
$X_i$ having simple polar locus along the (possibly reducible)
curve
$$C_i = X_i \cap \overline{(X \setminus X_i)} = \sum_{j \neq i} C_{ij},$$ for
each $1 \le i \le v$ (recall that $C_{ij} = C_{ji}$ and
that $C_i \neq \emptyset$ for each $ 1 \leq i \leq v$, because of
the connectedness hypothesis of $X$).

Take an index pair $i < j$ such that $C_{ij} \neq \emptyset$ and
consider $C_{ij}$, which is both an irreducible component of $C_i$
and of $C_j$. As in \eqref{eq:presnod}, take $\R_{C_{ij}}
(\omega_i)$ the Poincar\'e residue of $\omega_i$ on $C_{ij}$ and
denote it by $ \omega_{ij}$. Then, we have
\begin{equation}\label{eq:wij=-wjibis}
\omega_{ij}=-\omega_{ji},
\end{equation}for each $1 \leq i < j \leq v$.
In the trivial case $C_{ij} = \emptyset$, we have
$\omega_{ij}=\omega_{ji} = 0$; so \eqref{eq:wij=-wjibis} holds.

Fix an index pair $i<j$ such that $C_{ij} \neq \emptyset$. By
recalling our Remark \ref{rem:indices} and by our hypotheses, when
$X_i \cap X_j \cap X_k \neq \emptyset$ the lexicographically
ordered index triple uniquely determine the triple point set
either $\Sigma_{kij}$, or $\Sigma_{ikj}$, or $\Sigma_{ijk}$ on
$C_{ij}$, according to the case that either $ k < i$, or $i < k <
j$, or $k > j$. Now, given $i <j$ as above, we can consider
$\omega_{ij}$ which is a meromorphic 1-form on the curve $C_{ij}
\subset X_i$ having simple poles at the points in $\Sigma_{kij},
\Sigma_{ikj}, \Sigma_{ijk} \subset C_{ij}$ defined above and
determined by those $k \neq i, j$ such that $X_i \cap X_j \cap X_k
\neq \emptyset$. Otherwise, when $k$ is such that $\Sigma_{ijk}=
\emptyset$, then $\omega_{ij}$ must be considered as holomorphic
at $\Sigma_{ijk}$ and its residues at the empty point set are
zero; such $k$'s are determined by those vertices $v_k$ in $G_X$
which do not form any $3$-face with the edge $e_{ij}= (v_i, v_j)$.

On the other hand, if the index pair $i <j$ is such that $C_{ij} =
\emptyset$, one has that $\Sigma_{kij} = \Sigma_{ikj} =
\Sigma_{ijk} = \emptyset$ for each $k \neq i, j$ and that each
$\omega_{ij}$ is zero.

Therefore, for simplicity of notation, for any index pair $i <j$
we write
$$\omega_{ij} \in H^0 \biggl(C_{ij}, \omega_{C_{ij}} \biggl(\sum_{k < i}
\Sigma_{kij} + \sum_{k \in (i,j) } \Sigma_{ikj} + \sum_{k > j}
\Sigma_{ijk} \biggr)\biggr),$$ recalling that, if $C_{ij} =
\emptyset$, the above is the zero-vector space otherwise, if
$C_{ij} \neq \emptyset$ but some triple point set is the empty
set, its points do not impose any pole to the meromorphic form
$\omega_{ij}$.

In any case, one can compute the Poincar\'e residues at the given
points, namely $\R_{P^r_{kij}} (\omega_{ij})$, $\R_{P^s_{ikj}}
(\omega_{ij})$ and $\R_{P^t_{ijk}} (\omega_{ij})$, for any $ 1
\leq r \leq m_{kij}$, $1 \leq s \leq m_{ikj}$ and $1 \leq t \leq
m_{ijk}$.

To simplify our notation, if e.g.\ $k > j$, we write
$(\omega^t_{ij})_k$ (or directly $\omega^t_{ijk}$) to denote the
Poincar\'e residue $\R_{P^t_{ijk}} (\omega_{ij})$ of $\omega_{ij}$
at the $t^{\rm th}$ point $P^t_{ijk}$ of the set $\Sigma_{ijk}$, for
any $ 1 \leq t \leq m_{ijk}$. Similar notation for the other two
cases.

As observed in Remark \ref{rem:indices}, if we fix $i<j<k$ we
focus on the triple point set $\Sigma_{ijk}$ of $X$, which is
given by
$$X_i \cap X_j \cap X_k = C_{ij} \cap C_{ik} \cap C_{jk},$$
where
$$C_{ij}, C_{ik} \subset X_i, \;\;\; C_{ji} = C_{ij} , C_{jk} \subset X_j,
\;\;\; C_{ki} = C_{ik} , C_{kj} = C_{jk} \subset X_k.$$Therefore,
at any given triple point $P^t_{ijk}$, with $i < j < k$ and $1
\leq t \leq m_{ijk}$, one can compute six different residues.
Indeed, once we choose one of the three surfaces as the ambient
variety, we have two different possible choices of smooth,
irreducible curves (and so of meromorphic $1$-forms) to use for
such a computation. By taking into account \eqref{eq:wij=-wjibis},
Remark \ref{rem:residuinodi} and the lexicographic order of the
indices, the residues at $P^t_{ijk}$ satisfy
\begin{equation}\label{eq:wijkbis}
\omega^t_{\sigma(i)\sigma(j)\sigma(k)} = \sgn(\sigma) \;
\omega^t_{ijk},
\end{equation}where $i<j<k$, $1 \leq t \leq m_{ijk}$ and where
$\sigma \in \Sym(\{i,j,k \})$. Therefore, for each section $\omega
\in H^0 (X, \omega_X)$ there is, up to sign, a well determined
value associated to each point $P^t_{ijk} \in \Sigma_{ijk}$.
Recall that each such value is zero either if $\Sigma_{ijk}$ is
the empty point set or if some of the double curves is the empty
set.

For $i < j$ such that $C_{ij} \neq \emptyset$, the subsets of
triple points of $X$ lying on the curve $C_{ij}$ are parametrized
by those indices $k \neq i,j$ such that the vertex $v_k$ of $G_X$
forms a number bigger than or equal to one of faces with the edge
$e_{ij}=(v_i,v_j)$. By the Residue theorem on $C_{ij}$ and by
\eqref{eq:wijkbis}, we get
\begin{equation}\label{eq:resicurvebis}
\sum_{k \neq i,j}\sum_{t
=1}^{m_{\sigma_k(i)\sigma_k(j)\sigma_k(k)}} \sgn(\sigma_k) \;
\omega^t_{\sigma_k(i)\sigma_k(j)\sigma_k(k)} = 0,
\end{equation}where $\sigma_k \in \Sym (\{i,j,k \})$ is the permutation which
lexicographically reorders the index triples $(i,j,k)$.

Otherwise, if $C_{ij} = \emptyset$, then \eqref{eq:resicurvebis}
is trivially true, i.e.\ a sum of zeroes equals zero.

Choose, once and for all, the lexicographic orientation on the
graph $G_X$. Therefore, if the edge $\lambda_{ij}=(v_i,v_j)$
belongs to $G_X$, then it is an arrow from $v_i$ to $v_j$ if and
only if $i<j$. Recall that a set $\Sigma_{ijk}$, for some $i<j<k$,
corresponds to $m_{ijk} = |\Sigma_{ijk}|$ faces of $G_X$ insisting
on the triple of vertices $\{ v_i, v_j, v_k \}$. For any $1 \leq t
\leq m_{ijk}$, we associate to the $t^{\rm th}$-face the residue
$\omega^t_{ijk}$ computed as above.

From \eqref{eq:wijkbis}, \eqref{eq:resicurvebis} and from
the fact that $G_X$ is a $2$-dimensional graph, it follows that
the above computations determine a 2-cycle $\{\omega^t_{ijk}\}$ of
the graph $G_X$.

To sum up, the map $f$ is defined as a composition of maps in the
following way:
\begin{align*}
H^0(\omega_X) & \stackrel{i}{\hookrightarrow} \oplus
H^0(\omega_{X_i}(C_i)) \mapright{a}
\oplus H^0(\omega_{C_i}) \mapright{b} \\
& \mapright{b} \oplus H^0(\omega_{C_{ij}}(\sum_{k,t} P^t_{ijk}))
\mapright{c} \oplus H^0(\omega_{C_{ij}}(\sum_{k,t}
P^t_{ijk}))_{|P^t_{ijk}}) \to H_2(G_X, \CC)
\end{align*}
where $i$ is the natural inclusion and $a$, $b$, $c$ are given by
the following exact sequences:
\begin{align*}
&0 \to \oplus H^0(\omega_{X_i}) \to \oplus H^0(\omega_{X_i}(C_i))
\mapright{a} \oplus H^0(\omega_{C_i}) \\
&0 \to \oplus H^0(\omega_{C_i}(-C_{ij})) \to \oplus
H^0(\omega_{C_i}) \mapright{b} \oplus
H^0(\omega_{C_i|C_{ij}})\cong\oplus
H^0(\omega_{C_{ij}}(\sum_{k,t} P^t_{ijk})) \\
&0 \to \oplus H^0(\omega_{C_{ij}}(\sum_{l \neq k, r \neq t}
P^r_{ijl})) \to \oplus H^0(\omega_{C_{ij}}(\sum_{l , r }
P^r_{ijl})) \mapright{c} \oplus H^0(\omega_{C_{ij}}(\sum_{l , r }
P^r_{ijl}))_{|P^t_{ijk}})
\end{align*}

To compute $\ker(f)$, take as above $C_i = X_i \cap
\overline{(X \setminus X_i)} = \sum_{j\neq i} C_{ij}$, where we
recall that $C_{ij} = C_{ji}$, for $i \neq j$, some --- but not
all --- of them possibly empty. As in formulas \eqref{eq:presnod}
and \eqref{eq:fid}, denote by $\R_{C_i}$ the Poincar\'e residue
map on the reducible, nodal curve $C_i \subset X_i$ and by
$\delta_i$ the coboundary map on the surface $X_i$, for $ 1 \leq i
\leq v$. Thus, the first row of the diagram:

\begin{equation}\label{eq:diagrambis}
\xymatrix@C+1mm{%
0 \ar[r]
 & \displaystyle \bigoplus_{i=1}^v H^0(\omega_{X_i}) \ar[r] \ar[dr]!UL_-{\zeta}
   & \displaystyle \bigoplus_{i=1}^v H^0(\omega_{X_i}(C_i)) \ar[r]^-
{\oplus \R_{C_i}}
     & \displaystyle \bigoplus_{i=1}^v H^0(\omega_{C_i}) \ar[r]^-
{\oplus \delta_{i}}
       & \displaystyle \bigoplus_{i=1}^v H^1(\omega_{X_i}) \\
  &  & \ker(f) \ar[u]^-{\iota} \ar[r]^-{\sigma_{1/2}}
       & \displaystyle \bigoplus_{i<j} H^0(\omega_{C_{ij}}) \ar[u]^-\beta
\ar[]!UR;[ur]!DL_-\Delta }
\end{equation}
is naturally defined and exact. Apart from $\Delta$, introduced in
\eqref{eq:4.fid}, our aim is to define the maps $\varsigma$,
$\iota$, $\sigma_{1/2}$ and $\beta$ in such a way that the whole
diagram commutes and that the subsequence determined by the maps
$\varsigma$, $\sigma_{1/2}$ and $\Delta$ is exact.

Obviously, $\varsigma$ and $\iota$ are natural inclusions by the
very definition of $H^0(X, \omega_X)$. For what concerns the map
$\beta$, it suffices to define its image in one direct summand,
i.e.
$$(\pi_h \circ \beta) : \bigoplus_{i<j} H^0(\omega_{C_{ij}}) \to H^0
(\omega_{C_h}),$$ where $\pi_h$ is the projection on the
$h^{\rm th}$-summand, for a given $h \in \{1, \ldots, v\}$. When $h
\neq i,j$, the image is $0$, therefore the relevant summands are
the following:
\begin{equation}\label{eq:bbis}
\begin{aligned}
\bigoplus_{h<j} H^0(\omega_{C_{hj}}) \oplus \bigoplus_{i<h}
H^0(\omega_{C_{ih}})
 &  \longrightarrow H^0 (\omega_{C_h}), \\
\biggl(\bigoplus_{h<j} \omega_{hj}\;\;\;\; ,\;\;\;\;  \bigoplus_{i<h} \omega_{ih}
\biggr)
  & \mapsto \sum_{h<j} \omega_{hj}- \sum_{i<h} \omega_{ih},
\end{aligned}
\end{equation}
with the obvious condition that $\omega_{lm} = 0$ when $C_{lm} =
\emptyset$. First of all observe that $\beta$ is well-defined by
the definition of $H^0 (C_h, \omega_{C_h})$ (see Remark
\ref{rem:residuinodi}); moreover, the coefficients $\pm 1$ are
uniquely determined by the fact that $C_h \subset X_h$ and by
Remark \ref{rem:segnocurve}.

To define $\sigma_{1/2}$, recall that $\ker (f) \subseteq
H^0(X,\omega_X) \subset \oplus_{i=1}^v H^0 (\omega_{X_i} (C_i))$;
thus, an element in $\ker(f)$ is a collection of $v$ meromorphic
$2$-forms $(\gamma_1, \ldots, \gamma_v) \in \oplus_{i=1}^v H^0
(\omega_{X_i} (C_i))$ such that
$$\gamma_{ij} = - \gamma_{ji}, \; {\rm for} \; i < j,$$and
$$\gamma^t_{\tau_k(i) \tau_k(j)\tau_k(k)} = 0 \; {\rm for \; each} \; k \neq i,
\; j \; {\rm and \; for \; each} \; 1 \leq t \leq m_{\tau_k(i)
\tau_k(j)\tau_k(k)},$$where $\tau_k \in \Sym(\{i,j,k \})$ is the
permutation which lexicographically reorders the index triple
$\{i,j,k \}$. As before, we can limit ourselves to define its
image on a given direct summand; therefore, if $\pi_{ij}$ is the
projection on the $(ij)^{\rm th}$-summand, with $i < j$, then we have
the following equivalent expressions
\begin{equation}\label{eq:sigmabis}
(\pi_{ij} \circ \sigma_{1/2}) (\gamma_1, \ldots, \gamma_v) =  \frac{1}{2}
(\R_{C_{ij}}(\gamma_i) - \R_{C_{ij}}(\gamma_j)) =
\R_{C_{ij}}(\gamma_i) = - \R_{C_{ij}}(\gamma_j) .
\end{equation}
Observe that $\sigma_{1/2}$ is well-defined
since the $\gamma_i$'s are in the kernel of $f$; furthermore, when
$C_{ij} = \emptyset$, the image is obviously $0$.

By using the definition of $\Delta$ as in \eqref{eq:fid2}, it is
straightforward to check that diagram \eqref{eq:diagrambis}
commutes. Furthermore, it is trivial to show that
$$\im(\varsigma) = \ker(\sigma_{1/2})  \text{ and }
\im(\sigma_{1/2}) \subseteq \ker(\Delta).$$ To show the converse,
take $\alpha \in \ker (\Delta)$, thus $(\oplus_i
\delta_i(\beta(\alpha))) = 0$, i.e.\ $\beta(\alpha) \in \im(
\oplus_i (\R_{C_i}))$. This implies that $\alpha \in
\im(\sigma_{1/2}) $.

From the fact that the subsequence in \eqref{eq:diagrambis} is exact, it
follows that
$$\ker(f) \cong \ker( \sigma_{1/2}) \oplus \im(\sigma_{1/2}) \cong
\im(\varsigma) \oplus \ker(\Delta) \cong \bigoplus_{i=1}^d
H^0(\omega_{X_i})\oplus\ker(\Delta).$$ By Proposition
\ref{prop:Fi}, it follows that
\begin{equation}\label{eq:kerf3bis}
\ker(f) \cong \bigoplus_{i=1}^d H^0(\omega_{X_i}) \oplus
\coker(\Phi).
\end{equation}
This proves \eqref{eq:Fpgbis}.

To show that the map $f$ is surjective, we have to reconstruct a
global section of $\omega_X$ once we have a collection $ \{
\omega^t_{ijk} \} \in H_2(G_X, \CC)$.

Fix two indices $l < m $ in $\{1, \ldots, d \}$, such that $C_{lm}
\neq \emptyset$; this means that we consider the curve $C_{lm}$ as
the ambient variety to make our computations. From \eqref{eq:wijkbis}, we have
three different possibilities:

\begin{itemize}
\item if $k < l$, then $ \R_{P^t_{klm}}(\omega_{lm}) =
\omega^t_{lmk} = \sgn((k,l,m)) \omega^t_{klm} = \omega^t_{klm}$
for any $P^t_{klm} \in \Sigma_{klm}$, where $(k,l,m)$ is a
$3$-cycle in $\Sym(\{k,l,m\})$; \item if $l< k < m$, then $
\R_{P^t_{lkm}}(\omega_{lm})= \omega^t_{lmk} = \sgn((m,k))
\omega^t_{lkm} = - \omega^t_{lkm}$ for any $P^t_{lkm} \in
\Sigma_{lkm}$, where $(m,k)$ is a transposition in
$\Sym(\{l,k,m\})$; \item if $k > m$, we directly have
$\R_{P^t_{klm}}(\omega_{lm}) =\omega^t_{lmk} $ for any $P^t_{lmk}
\in \Sigma_{lmk}$.
\end{itemize}Therefore, by \eqref{eq:resicurvebis}, on $C_{lm}$ we have:
$$\sum_{k<l} \sum_{t=1}^{m_{klm}}\omega^t_{klm}  - \sum_{k \in (l,m) \cap \NN}
\sum_{t=1}^{m_{lkm}}\omega^t_{lkm} + \sum_{k>m}
\sum_{t=1}^{m_{lmk}}\omega^t_{lmk} =0.$$ This means that the
divisor
\begin{equation}\label{eq:omologobis}
D = \sum_{k<l} \sum_{t=1}^{m_{klm}}\omega^t_{klm} P^t_{klm}
  - \sum_{k \in (l,m) \cap \NN} \sum_{t=1}^{m_{lkm}}\omega^t_{lkm} P^t_{lkm}
  + \sum_{k>m} \sum_{t=1}^{m_{lmk}}\omega^t_{lmk} P^t_{lmk} \in Div (C_{lm})
\end{equation}is homologous to zero. By the Residue Theorem,
\eqref{eq:omologobis} implies there exists a global meromorphic
$1$-form $\omega_{lm} \in H^0 (C_{lm}, \omega_{C_{lm}} (D))$
having the given residues at the points in $Supp(D)$, i.e.\ such
that
$$\R_{P^r_{klm}}(\omega_{lm}) = \omega^r_{klm}, \;
\R_{P^s_{lkm}}(\omega_{lm}) = - \omega^s_{lkm}, \;
\R_{P^t_{lmk}}(\omega_{lm}) = \omega^t_{lmk},$$where $1\leq r \leq
m_{klm}$, $1\leq s \leq m_{lkm}$ and $1\leq t \leq m_{lmk}$.

The above discussion obviously holds for each choice of index
pairs. Fix now an index $1 \leq h \leq v$ and consider on the
surface $X_h$ the reducible nodal curve $C_h = X_h \cap
\overline{(X \setminus X_h)}$; since we are on $X_h$, by Remark
\ref{rem:segnocurve}, we can write $$C_h = C_h^+  + C_h^-,$$where
$$C_h^- = \sum_{l < h} C_{lh}, \; {\rm and} \;  C_h^+ = \sum_{ m > h}
C_{hm}.$$Thus, by the above discussion, on each $C_{lh}$
(resp.\ $C_{hm}$) we have a meromorphic $1$-form $ - \omega_{lh}$
(resp.\ $\omega_{hm}$) inducing the given residues at the given triple
points. Fix three indeces $i, j, h$ and consider the set of triple
points given by $X_i \cap X_j \cap X_h \neq \emptyset$. Since we
are on the surface $X_h$, such set of points is determined by the
intersection of $C_{ih} = C_{hi}$ and $C_{jh} = C_{hj}$. We have
the following possibilities:

\begin{itemize}
\item if $i < j < h$, the set of triple points is $\Sigma_{ijh}$
and on $C_{ih}$ (resp., on $C_{jh}$) we have the meromorphic
$1$-form $- \omega_{ih}$ (resp., $- \omega_{jh}$); therefore, by
\eqref{eq:wijkbis},
$$ \R_{P^t_{ijh}}(- \omega_{ih}) + \R_{P^t_{ijh}}(- \omega_{jh}) =
\omega^t_{ijh} - \omega^t_{ijh}= 0,$$for any $1 \leq t \leq
m_{ijh}$; \item if $h < i < j $, the set of triple points is
$\Sigma_{hij}$ and on $C_{hi}$ (resp., on $C_{hj}$) we have the
meromorphic $1$-form $\omega_{hi}$ (resp., $\omega_{hj}$); as
before,
$$ \R_{P^t_{hij}}(\omega_{hi}) + \R_{P^t_{hij}}(\omega_{hj}) =
\omega^t_{hij} - \omega^t_{hij}= 0,$$for any $1 \leq t \leq
m_{hij}$; \item if $i < h < j $, the set of triple points is
$\Sigma_{ihj}$ and on $C_{ih}$ (resp., on $C_{hj}$) we have the
meromorphic $1$-form $- \omega_{ih}$ (resp., $\omega_{hj})$; thus,
$$ \R_{P^t_{ihj}}(- \omega_{ih}) + \R_{P^t_{ihj}}(\omega_{jh}) =
- \omega_{ihj} + \omega_{ihj}= 0,$$for any $1 \leq t \leq
m_{ihj}$.
\end{itemize}

In either case, by Remark \ref{rem:residuinodi}, such forms glue
together to determine an element in $H^0 (C_{h}, \omega_{C_{h}})$.
This can be done for each $1 \leq h \leq d$, determining a
collection $\{ \overline{\omega}_h \} \in \bigoplus_{h=1}^v H^0
(\omega_{C_h})$.

Assume now to be in the case of hypothesis $(i)$, so each $X_h$ is
regular; since $h^1(X_h, \omega_{X_h}) = 0$, for each $1 \leq h
\leq v$, by the exact sequences
$$ 0 \to \omega_{X_h} \to \omega_{X_h} (C_h) \to \omega_{C_h} \to 0, \; 1 \leq h
\leq d,$$the collection of forms $\{ \overline{\omega}_h \}$ lifts
up to a collection of forms  $\{ \omega_h \} \in \bigoplus_{h=1}^v
H^0 (\omega_{X_h}(C_h))$. Take an index pair with $h < k$; since
$C_{hk}$ is both a component of $C_h^+$ on $X_h$ and of $C_k^-$ on
$X_k$, then
$$ \R_{C_{hk}}(\omega_k) = - \R_{C_{hk}}(\omega_h).$$This means that the
collection $\{ \omega_h \}$ is an element of $H^0(X, \omega_X)$,
so the map $f$ is surjective and formula \eqref{eq:Fpgbis} is
proved.

Assume now to be in the case of hypothesis $(ii)$; let $X_{j}$ be
an irreducible component of $X$, which is assumed to be an
irregular surface. By Dolbeault cohomology, by the hypothesis on
$C_{j}$ and by the Kodaira vanishing theorem, we get the following
diagram:
\begin{equation}\label{eq:diagram2bis}
\xymatrix{%
H^0(\omega_{X_{j}}(C_{j})) \ar[r]^-{\R_{C_j}} &
H^0(\omega_{C_{j}})
\ar[r]^{\delta_j} & H^1(\omega_{X_{j}}) \ar[r] & 0 \\
 & & H^0(X_{j}, \Omega^1_{X_{j}}) \ar[]!UL;[ul]!DR^-{\tr_{C_j}} \ar[u]^-{\eta_j}
},
\end{equation}which can be seen to be commutative,
where $\tr_{C_{j}}$ is the trace map of holomorphic $1$-forms on
$X_{j}$ to holomorphic $1$-forms on $C_{j}$ whereas $\eta_{j}$ is
the map defined by the cup product with the first Chern class of
$C_{j}$, $\Upsilon_{C_{j}} \in H^{1,1} (X_{j})$.

From the
surjectivity of $\delta_{j}$, it follows that
$$H^0( \omega_{C_{j}}) \cong \im(\R_{C_{j}}) \oplus
H^1(\omega_{X_{j}});$$on the other hand, by the Hard Lefschetz
theorem, $\eta_{j}$ is an isomorphism. Since $\delta_{j}$ is
injective on $\im(\tr_{C_{j}})$, then $\im(\R_{C_{j}}) \cap
\im(\tr_{C_{j}}) = \{ 0 \}$. On the other hand,
$$\im(\tr_{C_{j}}) \to H^1(\omega_{X_{j}})$$
is surjective. Then,
$$H^0( \omega_{C_{j}}) \cong \im(\R_{C_{j}}) \oplus \im(\tr_{C_{j}}).$$
Observe that the elements of $\im(\tr_{C_{j}})$ give zero residues
at the triple points of $X$ lying on $X_{j}$, since such
elements are restrictions to $C_{j}$ of global holomorphic 1-forms
of $X_{j}$. Therefore, the element $\omega_{j} \in
H^0(\omega_{C_{j}})$ of the collection $\{ \overline{\omega}_h \}
\in \bigoplus_{h=1}^v H^0 (\omega_{C_h})$, which was constructed
from the given non-zero collection of residues in $H_2(G_X)$, is
determined via $\R_{C_{j}}$ by an element in $H^0
(\omega_{X_{j}}(C_{j}))$, which is necessarily not zero. Then we
can conclude as above, proving also in this case the surjectivity
of $f$.
\end{proof}

In case $X$ is a planar Zappatic surface,
Theorem \ref{thm:4.pgbis} implies the following:

\begin{coroll}\label{cor:pgplanesbis}
Let $X = \bigcup_{i=1}^v \Pi_i$ be a planar Zappatic surface which
has only $E_3$ points as Zappatic singularities. Then,
\begin{align}
\label{eq:pgplanesbis}
 p_g(X) & = b_2(G_X), \\
\label{eq:pgplanesbis2}
 q(X) & = b_1(G_X).
\end{align}
\end{coroll}
\begin{proof}
Formula \eqref{eq:pgplanesbis} trivially follows from Theorem \ref{thm:4.pgbis}.
Notice that, in such a case, the proof of Theorem \ref{prop:4.chi} becomes simpler.
Indeed, each $\Sigma_{ijk}$ (cf.\ Notation \ref{not:Cij}) is either a singleton or empty,
since the double curves are lines (cf.\ Remark \ref{rem:planes}).

Formula \eqref{eq:pgplanesbis2} follows from \eqref{eq:chiplan},
\eqref{eq:pgplanesbis} and from the fact that
$\chi(G_X) = 1 - b_1(G_X) + b_2(G_X)$.
\end{proof}

In \S \ref{S:zapdeg} we shall extend the above results to
a good planar Zappatic surface $X$, assuming that $X$ is smoothable,
i.e.\ the central fibre of a degeneration.

\section{Zappatic degenerations}\label{S:zapdeg}

In this section we will focus on degenerations of smooth surfaces
to Zappatic ones.

\begin{definition}\label{def:degen}
Let $\D$ be the spectrum of a DVR (equiv.\ the complex unit disk).
Then, a {\em degeneration} (of relative dimension $n$) is a proper
and flat algebraic morphism
\[
\xymatrix@R=5mm{\X \ar[d]^\pi \\
\Delta}
\]
such that $\X_t = \pi^{-1}(t)$ is a smooth, irreducible,
$n$-dimensional projective variety, for $t \neq 0$.

If $Y$ is a smooth, projective variety, the degeneration
\[
\xymatrix@C=0mm@R=5mm{\X \ar[d]^\pi & \subseteq & \Delta \times Y \\
\Delta}
\]
is said to be an {\em embedded degeneration} in $Y$ of relative
dimension $n$. When it is  clear from the context, we will omit
the term embedded.

A degeneration (equiv.\ an embedded degeneration) is said to be
{\em semistable} if the total space $\X$ is smooth and if the
central fibre $\X_0$ (where $0$ is the closed point of $\D$) is a
divisor in $\X$ with global normal crossings, i.e.\ $\X_0 = \sum
V_i$ is a sum of smooth, irreducible components $V_i$'s which meet
transversally so that locally analitically the morphism $\pi$ is
defined by
\[
(x_1, \ldots, x_{n+1})
\xrightarrow{\;\pi\;} x_1x_2 \cdots x_k =t \in \D, \;\; k \leq
n+1.
\]
\end{definition}\medskip

Given an arbitrary degeneration $\pi:\X\to\D$, the well-known
Semistable Reduction Theorem (see \cite{Kempf}) states that there exists
a base change $\beta : \D \to \D$ (defined by $\beta(t) = t^m$,
for some $m$), a semistable degeneration $\psi : {\mathcal Z} \to
\D$ and a diagram
\[
\xymatrix@C-1mm@R-1mm{\Z \ar@{-->}[r]^f \ar[dr]_\psi & \X_\beta
\ar[d] \ar[r] &
\X \ar[d] \\
& \Delta \ar[r]^\beta & \Delta }
\]
such that $f$ is a birational map obtained by blowing-up and
blowing-down subvarieties of the central fibre.
Therefore, statements about degenerations which are invariant
under blowing-ups, blowing-downs and base-changes can be proved by
directly considering the special case of semistable degenerations.

From now on, we will be concerned with degenerations of relative
dimension two, namely degenerations of smooth, projective surfaces.

\begin{definition}\label{def:zappdeg}
Let $\X \to \D$ be a degeneration (equiv.\ an embedded
degeneration) of surfaces. Denote by $\X_t$ the general fibre,
which is by definition a smooth, irreducible and projective
surface; let $X=\X_0$ denote the central fibre. We will say that the
degeneration is {\em Zappatic} if $X$ is a Zappatic surface
and $\X$ is smooth except for:
\begin{itemize}
\item ordinary double points at points of the double locus of $X$,
which are not the Zappatic singularities of $X$;
\item further singular points at the Zappatic singularities of $X$ of type
$T_n$, for $n \geq 3$, and $Z_n$, for $n \geq 4$.
\end{itemize}
A Zappatic degeneration will be called \emph{good}
if the central fibre is moreover a good Zappatic surface.
Similarly, an embedded degeneration will be called a \emph{planar Zappatic degeneration}
if its central fibre is a planar Zappatic surface.

Notice that we require the total space $\X$ to be smooth
at $E_3$-points of $X$.
\end{definition}

If $\X\to \D$ is a good Zappatic degeneration, the singularities that
$\X$ has at the Zappatic singularities of the central fibre $X$
are explicitly described in \cite{CCFM}.

\begin{notation}\label{not:zappdeg}
Let $\X \to \D$ be a degeneration of surfaces and let $\X_t$ be
the general fibre, which is by definition a smooth, irreducible
and projective surface. Then, we consider some of the intrinsic
invariants of $\X_t$:
\begin{itemize}
\item $\chi : = \chi(\Oc_{\X_t})$;
\item $K^2 := K_{\X_t}^2$;
\item $p_g := p_g(\X_t)$;
\item $\chi_{\rm top} := \chi_{\rm top} (\X_t)$;
\end{itemize}If the degeneration is assumed to be embedded in $\PR$, for some
$r$, then we also have:
\begin{itemize}
\item $d : = \deg(\X_t)$;
\item $g := (K+H)H/2 + 1,$ the sectional genus of $\X_t$.
\end{itemize}
\end{notation}

We will be mainly interested in computing these invariants in
terms of the central fibre $X$. For some of them, this is quite simple.
For instance, when $\X \to \D$ is an embedded degeneration
in $ \PR$, for some $r$,  and if the central fibre $\X_0 = X =
\bigcup_{i=1}^v X_i$, where the $X_i$'s are smooth, irreducible
surfaces of degree $d_i$, $ 1 \leq i \leq v$, then
by the flatness of the family we have
$$
d = \sum_{i=1}^v d_i.
$$

When $\X \to \D$ is a good Zappatic degeneration (in particular a
good, planar Zappatic degeneration), we can easily compute some
of the above invariants by using our results of \S \ref{S:3}.
Indeed, by using our Notation \ref{def:vefg} and Propositions
\ref{prop:ghypsect}, \ref{prop:4.chi}, we get the following results.

\begin{proposition}\label{prop:gzapdeg}
Let $\X \to \D$ be a good Zappatic degeneration embedded in $
\PR$. Let $\X_0 = X = \bigcup_{i=1}^v X_i \subset \PR$ be the
central fibre and let $G = G_X$ be its associated graph. Let $C$
be the double locus of $X$, i.e.\ the union of the double curves of
$X$, $C_{ij} = C_{ji} = X_i \cap X_j$ and let $c_{ij} =
\deg(C_{ij})$. Let $D$ be a general hyperplane section of $X$
and let $D_i$ be the $i^{th}$ smooth, irreducible component of
$D$, which is a general hyperplane section of $X_i$, and denote
by $g_i$ its genus. Then:
\begin{align}
g &= \sum_{i=1}^v g_i + \sum_{e_{ij}\in E} \; c_{ij} -v +1.
\label{eq:gzapdeg}
\end{align}
When $X$ is a good, planar Zappatic surface, if
$G^{(1)}$ denotes the $1$-skeleton of $G$, then:
\begin{align}
g &= 1 - \chi(G^{(1)}) = e- v + 1 .\label{eq:gplanzapdeg}
\end{align}
\end{proposition}
\begin{proof}
It directly follows from our computations in Proposition
\ref{prop:ghypsect} and from the flatness of the family of
hyperplane sectional curves of the degeneration (cf.\ formula
\eqref{eq:genuscurves}).
\end{proof}

\begin{proposition}\label{prop:chizapdeg}
Let $\X \to \D$ be a good Zappatic degeneration and let $\X_0 = X
= \bigcup_{i=1}^v X_i $ be its central fibre. Let $G = G_X$ be its
associated graph and let $E$ the indexed set of edges of $G$. Let
$C$ be the double locus of $X$, which is the union of the double
curves $C_{ij} = C_{ji} = X_i \cap X_j$. Denote by $g_{ij}$ the
genus of the smooth curve $C_{ij}$. Then
\begin{align}
\chi &= \sum_{i=1}^v \chi(\Oc_{X_i}) - \sum_{e_{ij} \in E}
\chi(\Oc_{C_{ij}}) + f. \label{eq:chizapdeg}
\end{align}

Moreover, if $\X \to \D$ is a good, planar Zappatic
degeneration, then
\begin{align}
\chi &= \chi(G) = v -e + f. \label{eq:chiplanzapdeg}
\end{align}
\end{proposition}

\begin{proof}
It follows from Proposition \ref{prop:4.chi} and from the invariance of $\chi$ under flat degeneration.
\end{proof}

In the particular case that $\X \to \D$ is a semistable Zappatic degeneration,
i.e.\ if $X$ has only $E_3$-points as Zappatic singularities,
then $\chi$ can be computed also in a different way by topological methods
(see formula \eqref{eq:chiCS} in Theorem \ref{thm:4.CSbis}).

The above results are indeed more general:
$X$ is allowed to have any good Zappatic singularity,
namely $R_n$-, $S_n$- and $E_n$-points, for any $n\ge3$,
and moreover our computations do not depend on the fact that $X$ is smoothable,
i.e.\ that $X$ is the central fibre of a degeneration.
Notice also that a good Zappatic degeneration is not semistable in general.

For what concerns the geometric genus, assume now --- unless otherwise stated --- 
that the Zappatic surface $X=\bigcup_{i=1}^v X_i$ is
the central fibre of a semistable Zappatic degeneration
$\X\to\D$, i.e. $\X$ is smooth and $X$ has only $E_3$-points as Zappatic
singularities. In this case, Theorem \ref{thm:4.pgbis} implies the following:

\begin{proposition}\label{cor:.fibis}
Let $\X \to \D$ be a semistable Zappatic degeneration
and $\X_0=X=\bigcup_{i=1}^v X_i$ be its central fibre.
Let $G_X$ be the associated graph to $X$
and $\Phi$ be the map defined in \eqref{eq:4.fi}.
Then, for any $t\in\D$,
\begin{equation}\label{eq:minougpgbis}
p_g(\X_t) \leq b_2(G_X)+\sum_{i=1}^v p_g(X_i)+ \dim(\coker(\Phi)).
\end{equation}
\end{proposition}

\begin{proof}
By semi-continuity, we have $p_g(\X_t) \leq p_g(\X_0) = p_g(X)$. One concludes 
by using formula \eqref{eq:Fpgbis}.
\end{proof}

On the other hand, $p_g(\X_t)$ and $\chi(\X_t)$
can be also computed by topological methods:
indeed the Clemens-Schmid exact sequence relates the mixed Hodge theory
of the central fibre $X$ to that of $\X_t$ by means of the
monodromy of the total space $\X$
(see \cite{Morr} for definitions and statements).
In our particular situation, the following result holds:

\begin{theorem}[Clemens-Schmid]\label{thm:4.CSbis}
Let $\X \to \D$ be a semistable Zappatic degeneration
and $\X_0=X=\bigcup_{i=1}^v X_i$ be its central fibre.
Let $G_X$ be the associated graph to $X$
and $\Phi$ be the map defined in \eqref{eq:4.fi}.
Then, for any $t\ne0$,
\begin{align}
\label{eq:chiCS}
\chi(\Oc_{\X_t}) & = \chi(G_X), \\
\label{eq:pgCS}
p_g(\X_t)  & =   b_2(G_X) + \smash{\sum_{i=1}^d p_g(X_i)} + \dim(\coker (\Phi)).
\end{align}
\end{theorem}

A proof of Theorem \ref{thm:4.CSbis} can be found
in \cite[``Clemens-Schmid I and II"]{Morr}.
The above result, together with
Theorem \ref{thm:4.pgbis}, implies the following:

\begin{corollary}\label{cor:4.CSganzo}
Let $X=\bigcup_{i=1}^v X_i$ be the central fibre
of a semistable Zappatic degeneration $\X \to \D$.
Then, for every $t\in\D$,
\begin{equation}\label{eq:pgCS=}
p_g(\X_t) = p_g(X)= b_2(G_X) + \sum_{i=1}^v p_g(X_i) + \dim(\coker (\Phi)).
\end{equation}
In particular, the geometric genus of the fibres of $\X \to \D$ is \emph{constant}.
\end{corollary}

\begin{proof}
Formula \ref{eq:pgCS=} trivially follows from \eqref{eq:Fpgbis}, \eqref{eq:pgCS}
and from semicontinuity.
\end{proof}

Recalling the proof of Theorem \ref{thm:4.pgbis}, we also have the following:

\begin{corollary}\label{cor:4.CSganzo2}
Let $X=\bigcup_i X_i$ be a Zappatic surface with global normal crossings,
i.e.\ with only $E_3$-points as Zappatic singularities.
Let $G_X$ be its associated graph.
A necessary condition for the smoothability of $X$ is the surjectivity
of the homomorphism $f : H^0(X, \omega_X) \to H_2 (G_X)$
defined in \eqref{eq:kerfbis}.
\end{corollary}

\begin{proof}
If $X$ is smoothable, then \eqref{eq:pgCS=} implies that the equality
in \eqref{eq:Fpgbis} holds.
This implies that $f$ is surjective by the proof of Theorem \ref{thm:4.pgbis}.
\end{proof}

The above results naturally suggest the following:

\begin{question}
Is the homomorphism $f : H^0(X, \omega_X) \to H_2 (G_X)$
in \eqref{eq:kerfbis} always surjective? Equivalently,
does the equality in \eqref{eq:Fpgbis} always hold?
\end{question}

We believe that the answer to this question should be negative,
but we have not been able to exhibit a counterexample so far.

In case the answer to the above question were negative,
it should be interesting to compare the surjectivity of $f$
with other smoothability conditions, like Friedman's one in \cite{Frie}.

We now generalize the computations for $p_g$ to the case
of good Zappatic degenerations, i.e.\ degenerations where
the central fibre $X$ is a union of surfaces having not only
$E_3$-points, but also $R_n$-, $S_n$- and $E_n$-points for any $n\ge3$.

\begin{definition}\label{def:pg}
Let $\X\to\D$ be a good Zappatic degeneration and $X=\X_0$ be its central fibre.
Consider the semistable reduction $\X'\to\D$ of $\X\to\D$ together 
with its central fibre $X'=\X'_0$, which is a Zappatic surface with global normal crossings,
i.e.\ with only $E_3$-points.
We define the geometric genus as (cf.\ Remark \ref{rem:pg}):
\begin{equation}\label{eq:pgg}
p_g(X):=p_g(X')=h^0(X',\omega_{X'}).
\end{equation}
\end{definition}As we will see in a moment, the definition is well-posed.

\begin{theorem}\label{thm:pg}
Let $X=\X_0$ be a good Zappatic surface which is the central fibre
of a degeneration $\X\to\D$ and let $G_X$ be its associated graph.
Then
\begin{equation}\label{eq:pgCS=bis}
 p_g(X)= p_g(\X_t) = b_2(G_X) + \sum_{i=1}^v p_g(X_i) + \dim(\coker (\Phi)).
\end{equation}
\end{theorem}

\begin{proof}[Sketch of the proof]
Complete details will appear in \cite{CCFM}. Here we give an outline of the proof.

The first step is to understand how to get the semistable reduction locally
near $R_n$-, $S_n$- and $E_m$-points, for $n\ge3$ and $m\ge4$.

Consider a $R_n$-point $x\in X$.
Then $x$ is an isolated singularity for the total space $\X$
and it is a \emph{minimal singularity} in the sense of Koll\'{a}r (\cite{Kollar} and \cite{KolSB}).
Let $\tilde\X\to\X$ be the blow-up at $x$.
Its exceptional divisor $F$ is a minimal degree surface (of degree $n$)
in $\Pp^{n+1}=\Pp(\T_{\X,x})$, where $\T_{\X,x}$ is the tangent space of $\X$ at $x$.
Furthermore $F$ is connected in codimension one and can be explicitly described.
In particular one can show that either $F$ is smooth or $F$ has $R_m$-points, for $m\le n$. 
Some points $E_4$ can appear along the intersection of the exceptional divisor with the 
strict transform of $X$. In any event, after finitely many blow-ups of the total space at 
points, we resolve the singularities of the total space. It turns out that all the components of 
the exceptional divisor are rational and all the 
double curves involved in them are also rational.

The situation is completely similar for $S_n$-points and $E_n$-points, $n \geq 4$. 

Therefore one can get the semistable reduction $\X'\to\D$ of $\X\to\D$ just
by blowing-up $\X$ at points which are good Zappatic singularities of the central fibre.
Let $\sigma:\X'\to\X$ be this blow-up and 
$X'=\X'_0$ be the central fibre of $\X'\to\D$.

By \eqref{eq:pgg}, $p_g(X)$ is defined to be $p_g(X')$; 
Theorem \ref{thm:4.pgbis} tells how to compute it. 

Now, our second and last step is to prove that $p_g(X')$ is given by \eqref{eq:pgCS=bis} and 
that $p_g(X') = p_g(\X_t)$, for $t \neq 0$. Since the semistable reduction of $\X \to \D$ 
involves only the central fibre and since all the exceptional divisors as well as 
all the curves involved in them are 
rational, it suffices to prove that $b_2(G_X) = b_2 (G_{X'})$.

Consider an open $n$-face (resp.\ a $n$-angle) $G_x$ of $G_X$,
namely $G_x$ is the subgraph with $n$ vertices corresponding to $n$ planes
forming a $R_n$-point (resp.\ a $S_n$-point) $x$.
Let $G'_x$ be the subgraph of $G_{X'}$ containing the vertices corresponding
to the proper transforms of the $n$ planes
and the exceptional divisors contained in $\sigma^{-1}(x)$.
The above 
description of the infinitesimal neighbourhood of $x$ shows that, as topological spaces,
the subgraph $G_x$ is a deformation retract of $G_{X'}$.

Similarly, if $G_x$ is a closed $n$-face of $G_X$ (i.e.\ a subgraph with $n$ vertices corresponding
to $n$ planes forming an $E_n$-point $x$),
then the above description shows 
that the subgraph $G'_x$ of $G_{X'}$ containing the vertices corresponding to the proper transforms
of the planes and the exceptional divisors contained in $\sigma^{-1}(x)$ is a triangulation of $G_x$.

It follows that the graphs $G_X$ and $G_{X'}$ have the same homological invariants, which is what 
we had to prove.
\end{proof}

\end{document}